%% file: LinkingCombClassDyn.tex
\documentclass[12pt]{amsart}
\usepackage{amsmath}
\usepackage{amscd}
\usepackage{pb-diagram}
\usepackage{comment}
\usepackage{amssymb}
\usepackage{graphicx}
\usepackage{subcaption}

\input mathDefs

\newcommand{\mto}{\multimap}

\def\cell#1{\protect\mbox{$\stackrel{\circ}{#1}$}}
\def\ogr#1{\mbox{$\langle#1\rangle$}}
\def\csl#1{\mbox{$\langle#1\rangle_\lambda$}} 
\def\cse#1{\mbox{$\langle#1\rangle_\epsilon$}} 
\def\csb#1{\mbox{$\langle#1\rangle_\beta$}} 
\def\csd#1{\mbox{$\langle#1\rangle_\delta$}} 
\def\cstd#1{\mbox{$\langle#1\rangle_{\tilde{\delta}}$}} 
\def\csdp#1{\mbox{$\langle#1\rangle_{\delta'}$}} 
\newcommand{\sign}{\operatorname{sign}}

\newcommand{\Exit}{\operatorname{Ex}}  
\newcommand{\exit}{\operatorname{ex}}  
\newcommand{\ccl}{\operatorname{Cl}} 
\renewcommand{\Bd}{\Exit}
\newcommand{\ClasCon}{\operatorname{Con}}
\renewcommand{\omap}{\nrightarrow}
\renewcommand{\emptyset}{\varnothing}
\renewcommand{\subset}{\subseteq}

\newenvironment{aenum}{\begin{enumerate}

}{\end{enumerate}}

\renewcommand{\rho}{\varrho}
\renewcommand{\phi}{\varphi}
\renewcommand{\epsilon}{\varepsilon}
\newcommand{\currentDate}{\today}
\def\PP{\mathobj{\mathbb{P}}}
\newcommand{\isoS}{S(\cS)}

\articletheorems

\begin{document}

\author{Bogdan Batko}
\address{Bogdan Batko, Division of Computational Mathematics,
  Faculty of Mathematics and Computer Science,
  Jagiellonian University, ul.~St. \L{}ojasiewicza 6, 30-348~Krak\'ow, Poland
}
\email{Bogdan.Batko@uj.edu.pl}

\author{Tomasz Kaczynski}
\address{Tomasz Kaczynski, D\'epartement de math\'ematiques, Universit\'e de
Sherbrooke, 2500 boul. Universit\'e, Sherbrooke, Qc, Canada  J1K 2R1
}
\email{t.kaczynski@usherbrooke.ca}
\author{Marian Mrozek}
\address{Marian Mrozek, Division of Computational Mathematics,
  Faculty of Mathematics and Computer Science,
  Jagiellonian University, ul.~St. \L{}ojasiewicza 6, 30-348~Krak\'ow, Poland
}
\email{Marian.Mrozek@uj.edu.pl}
\author{Thomas Wanner}
\address{Thomas Wanner, Department of Mathematical Sciences,
George Mason University, Fairfax, VA 22030, USA
} \email{twanner@gmu.edu}
\date{today}
\thanks{Research of B.B.\ and M.M.\ was partially supported by
  the Polish National Science Center under Ma\-estro Grant No. 2014/14/A/ST1/00453.
  Research of T.K.\ was supported by a Discovery Grant from NSERC of Canada.
  T.W.\ was partially supported by NSF grants
  DMS-1114923 and DMS-1407087.
  All authors gratefully acknowledge the support of Hausdorff Research Institute for Mathematics in Bonn
  for providing an excellent environment to work together during the 2017 Special Hausdorff Program on
  Applied and Computational Algebraic Topology.}
\subjclass[2010]{Primary: 37B30; Secondary: 37E15, 57M99, 57Q05, 57Q15.}
 \keywords{Combinatorial vector field, multivalued dynamical system, simplicial complex,
 discrete Morse theory, Conley theory, Morse decomposition, Conley-Morse graph, isolated invariant set, isolating block.}

\title[Linking combinatorial and classical dynamics]
{Linking combinatorial and classical dynamics:
Conley index and Morse decompositions}

\date{Version compiled on \currentDate}

\begin{abstract}
We prove that every combinatorial dynamical system in the sense of Forman, defined on a family of simplices of a simplicial complex,
gives rise to a  multivalued dynamical system~$F$ on the geometric realization of the simplicial complex.
Moreover, $F$ may be chosen in such a way that the isolated invariant sets, Conley indices, Morse decompositions, and Conley-Morse graphs
of the two dynamical systems are in one-to-one correspondence.
\end{abstract}

\maketitle


\section{Introduction}
\label{sec:intro}

In the years since Forman \cite{Fo98a,Fo98} introduced combinatorial vector fields
on simplicial complexes, they have found numerous applications in
such areas as visualization and mesh compression \cite{LeLoTa04}, graph braid groups \cite{FaSa05},
homology computation \cite{HaMiMrNa2014,MiNa13},
astronomy \cite{So11}, the study of \v{C}ech and Delaunay complexes \cite{BaEd17},
and many others.
One reason for this success has its roots in Forman's original motivation. In his
papers, he sought to transfer the rich dynamical theories due to
Morse \cite{Mo34} and Conley \cite{conley:78a} from the continuous setting of a continuum (connected compact metric space)
to the finite, combinatorial setting of a simplicial complex.
This has proved to be extremely useful for establishing finite,
combinatorial results via ideas from dynamical systems.
In particular, Forman's theory yields an alternative when studying sampled dynamical systems.
The classical approach consists in the numerical study of the dynamics of the differential equation constructed
from the sample. The construction uses the data in the sample
either to discover the natural laws governing the dynamics \cite{SL07} in order to write the equations
or to interpolate or approximate directly the unknown right-hand-side of the equations \cite{BBCHS07}.
In the emerging alternative one can eliminate differential equations and study directly the combinatorial
dynamics defined by the sample \cite{Fo98a,Fo98,RoWoSh11,KRHHNDM16,Mr15}.

The two approaches are essentially distinct.
On the one hand, dynamical systems defined by differential equations
on a differentiable manifold
arise in a wide variety of applications and show an extreme wealth of
observable dynamical behavior, at the expense of fairly involved
mathematical techniques which are needed for their precise description.
On the other hand, the discrete simplicial complex setting makes the study of many
phenomena simple, due to the availability of fast combinatorial
algorithms. This leads to the natural question of which approach
should be chosen when for a given problem.

In order to answer this question it may be helpful to go beyond
the exchange of abstract underlying ideas present in much of
the existing work and look for the precise relation between the two theories.
In our previous paper \cite{KaMrWa16} we took this path and studied
the formal ties of multivalued dynamics in the combinatorial
and continuum settings. The choice of multivalued dynamics is natural,
because the combinatorial vector fields generate multivalued dynamics
in a natural way. Moreover, in the finite setting such dynamical phenomena as homoclinic
or heteroclinic connections are not possible in single-valued dynamics.
The choice of multivalued dynamics on continua is not a restriction.
This is a broadly studied and well understood theory.
The theory originated in the middle of the 20th century from the study of
contingent equations and  differential inclusions
 \cite{Wa61,Ro65b,AuCe1984} and control theory \cite{Ro65a}.
At the end of the 20th century it was successfully applied to computer assisted proofs in dynamics \cite{MiMr95,Mr96}.
In particular, the Conley theory for multivalued dynamics was studied by several authors
\cite{Mr90b,KaMr95,St06,DzKr08,DzGa11,BaMr16,Ba17}.

In \cite{KaMrWa16} we proved that for any combinatorial vector
field on the collection of simplices of a simplicial complex one can construct an acyclic-valued and
upper semicontinuous map on the underlying geometric realization
whose dynamics on the level of invariant sets exhibits the same complexity.
More precisely, by introducing
the notion of isolated invariant sets in the discrete setting, we
established a correspondence between isolated invariant sets in the
combinatorial and classical multivalued settings. We also presented a link
on the level of individual dynamical trajectories.

In the present paper we complete the program started in \cite{KaMrWa16}
by showing that the formal correspondence established there extends
to Conley indices of the corresponding isolated invariant set as well as
Morse decompositions and Conley-Morse graphs \cite{Ar09,BuAtAl2012},
a global descriptor of dynamics capturing its gradient structure.

The organization of the paper is as follows.
In Section~\ref{sec:main-result} we present the main result of the paper and
illustrate it with some examples. In Section~\ref{sec:conley} we recall the
basics of the Conley theory for multivalued dynamics.
In Section~\ref{sec:map-F} we recall from \cite{KaMrWa16} the construction
of a multivalued self-map $F:X \mto X$ associated with a combinatorial vector field $\cV$
on a simplicial complex $\cX$ with the geometric realization $X:=|\cX|$.
In Section~\ref{sec:correspondence} we use this construction to outline the proof of the main result of the paper
in a series of auxiliary theorems. The remaining sections are devoted to the proofs
of these theorems.

\section{Main result}
\label{sec:main-result}

Let $\cX$ denote the family of simplices of a finite abstract simplicial complex.
The face relation on $\cX$ defines on $\cX$ the  $T_0$ Alexandroff topology \cite{Al1937}.
A subset $\cA\subset \cX$ is {\em open} in this topology if
all cofaces of any element of $\cA$ are also in $\cA$.
The closure of $\cA$ in this topology, denoted $\ccl\cA$, is the family of all faces of all simplices in $\cA$
(see Section~\ref{sec:pre} for more details).
A {\em combinatorial vector field} $\cV$ on $\cX$ is a partition of $\cX$
into singletons and doubletons such that each doubleton consists of a simplex and
one of its cofaces of codimension one. The singletons are referred to as {\em critical cells}.
The doubleton considered as a pair with lower dimensional simplex coming first is referred to as a {\em vector}.
\begin{figure}[tb]
  \includegraphics[width=0.60\textwidth]{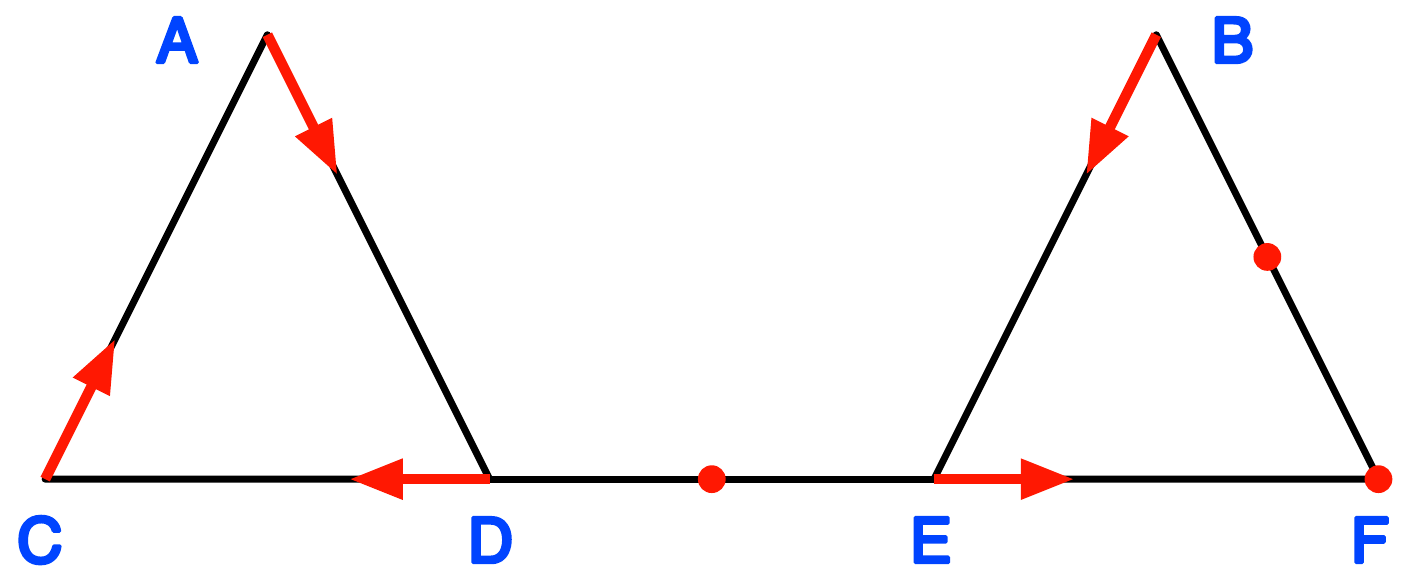}
  \caption{Sample discrete vector field. This figure shows a simplicial
           complex~$\cX$ which is a graph on six vertices with seven
           edges. Critical cells are indicated by red dots, vectors of
           the vector field are shown as red arrows.}
  \label{fig:graphexample}
\end{figure}

The elementary example in Figure~\ref{fig:graphexample} presents a one-dimensional
simplicial complex $\cX$ consisting
of six vertices $\{A,B,C,D,E,F\}$ and seven edges
$\{AC,AD,BE,BF,CD,DE,EF\}$, and the combinatorial vector field consisting of three singletons (critical cells)
$\{\{BF\},\{DE\},\{F\}\}$ and five doubletons (vectors) $\{\{A,AD\},\{B,BE\},\{C,AC\},\{D,CD\},\{E,EF\}\}$.
With a combinatorial vector field $\cV$ we associate multivalued dynamics
given as iterates of a multivalued map $\Pi_\cV:\cX\mto\cX$
sending each critical simplex to all of its faces, each source of a vector to the corresponding target, and each
target of a vector to all faces of the target other than the corresponding source and the target itself.
In the case of the
example in Figure~\ref{fig:graphexample} the map is
(we skip the braces in the case of singletons to keep the notation simple)
\begin{eqnarray*}
  \Pi_\cV&=&\{(A,AD),(AD,D),(B,BE),(BE,E),(BF,\{B,BF,F\}),\\
  &&~~(C,AC),(AC,A),(CD,C)(D,CD),(DE,\{D,DE,E\}),\\
  &&~~(E,EF),(EF,F),(F,F)\} \; .
\end{eqnarray*}
\begin{figure}[tb]
  \includegraphics[width=0.60\textwidth]{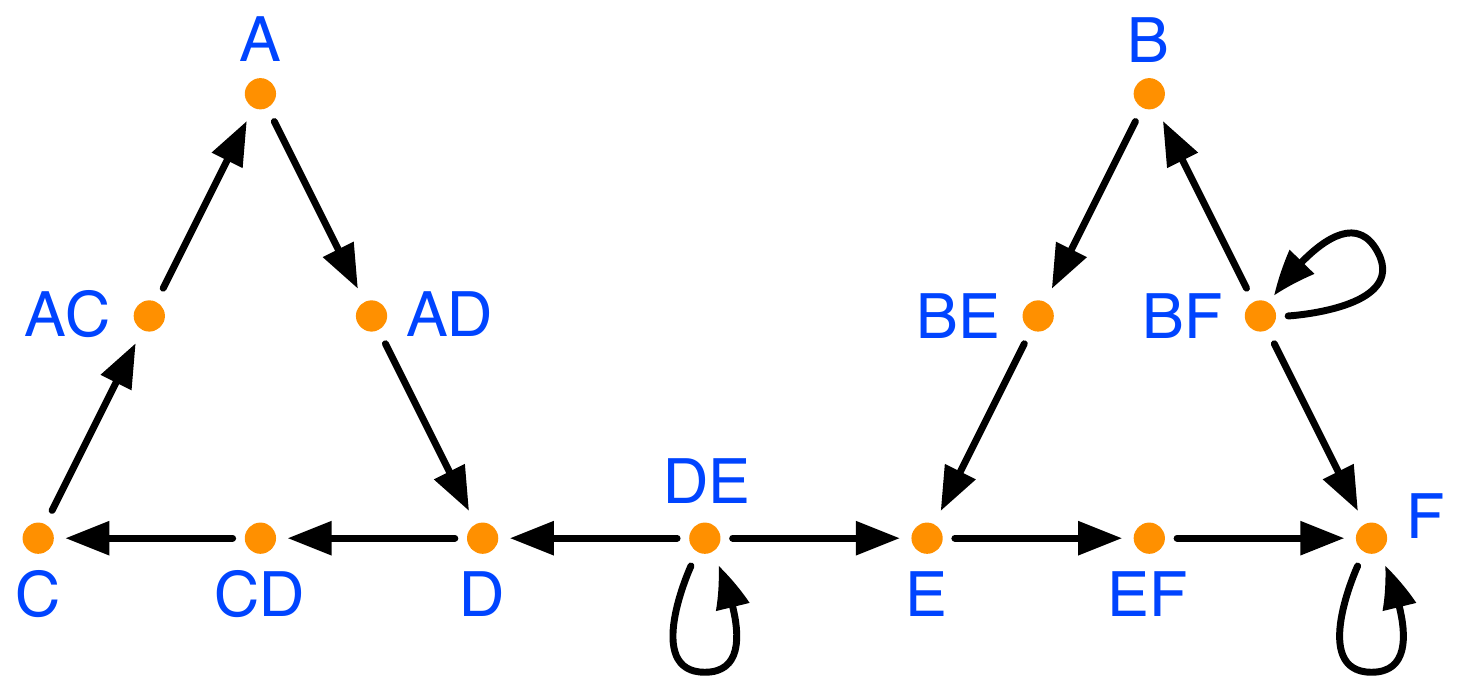}
  \caption{The directed graph $G_\cV$ for the combinatorial vector field in Figure~\ref{fig:graphexample}.}
  \label{fig:digraphexample}
\end{figure}

The multivalued map $\Pi_\cV$ may be considered as a directed graph $G_\cV$ with vertices in $\cX$ and an arrow
from a simplex $\sigma$ to a simplex $\tau$ whenever $\tau\in \Pi_\cV(\sigma)$.
The directed graph $G_\cV$ for the combinatorial vector field in Figure~\ref{fig:graphexample} is presented
in Figure~\ref{fig:digraphexample}.
A subset $\cA\subset\cX$ is {\em invariant} with respect to $\cV$ if every element of $\cA$ is both a head and a tail
of an arrow in $G_\cV$ which joins vertices in $\cA$.
An element $\sigma\in\ccl\cA\setminus\cA$ is an {\em internal tangency} of $\cA$ if it admits an arrow
originating in $\sigma$ with its head in $\cA$, as well as an arrow terminating in $\sigma$ with its tail in $\cA$.
The set $\Exit\cA:=\ccl\cA\setminus\cA$ is referred to
as the {\em exit set} of $\cA$ (see \cite[Definition 3.4]{KaMrWa16})
or {\em mouth} of $\cA$ (see \cite[Section 4.4]{Mr15}).
An invariant $\cS$ set is an {\em isolated invariant set} if the exit set $\Exit \cS$
is closed and it admits no internal tangencies.
Note that $\cX$ itself is an isolated invariant set if and only if it is invariant.
The {\em (co)homological Conley index} of an isolated invariant set $\cS$
is the relative singular (co)homology of the pair $(\ccl\cS,\Exit\cS)$.
Note that $(\ccl\cS,\Exit\cS)$ is a pair of simplicial subcomplexes
of the simplicial complex $\cX$. Therefore, by McCord's Theorem \cite{MC65},
the singular (co)homology of the pair $(\ccl\cS,\Exit\cS)$
isomorphic to the simplicial homology of the pair $(\ccl\cS,\Exit\cS)$.

The singleton $\{BF\}$ in Figure~\ref{fig:graphexample} is an example of an isolated invariant set of~$\cV$.
Its exit set is $\{B,F\}$ and its Conley index is the (co)homology of the pointed circle.
Another example is the set $\{A,AC,AD,C,CD,D\}$
with an empty exit set and the Conley index equal to the (co)homology of the circle.
Both these examples are minimal isolated invariant sets, that is, none of their proper non-empty subsets is
an isolated invariant set.

\begin{figure}[tb]
  \includegraphics[width=0.60\textwidth]{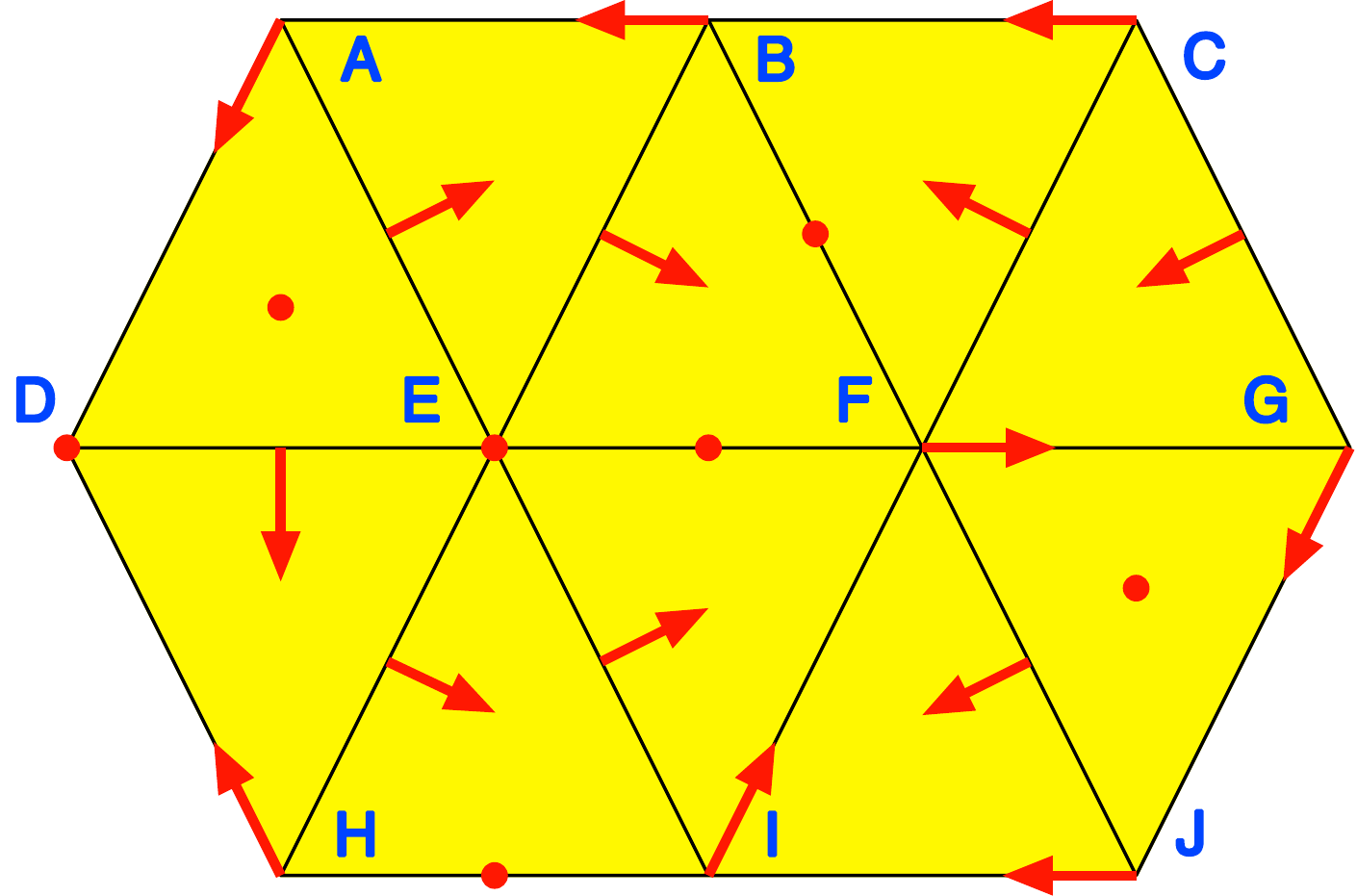}
  \caption{Sample discrete vector field. This figure shows a simplicial
           complex~$\cX$ which triangulates a hexagon (shown in yellow),
           together with a discrete vector field. Critical cells
           are indicated by red dots, vectors of the vector field
           are shown as red arrows. This example will be discussed
           throughout the paper.}
  \label{fig:mainexample}
\end{figure}

\begin{figure}[tb]
  \includegraphics[width=0.60\textwidth]{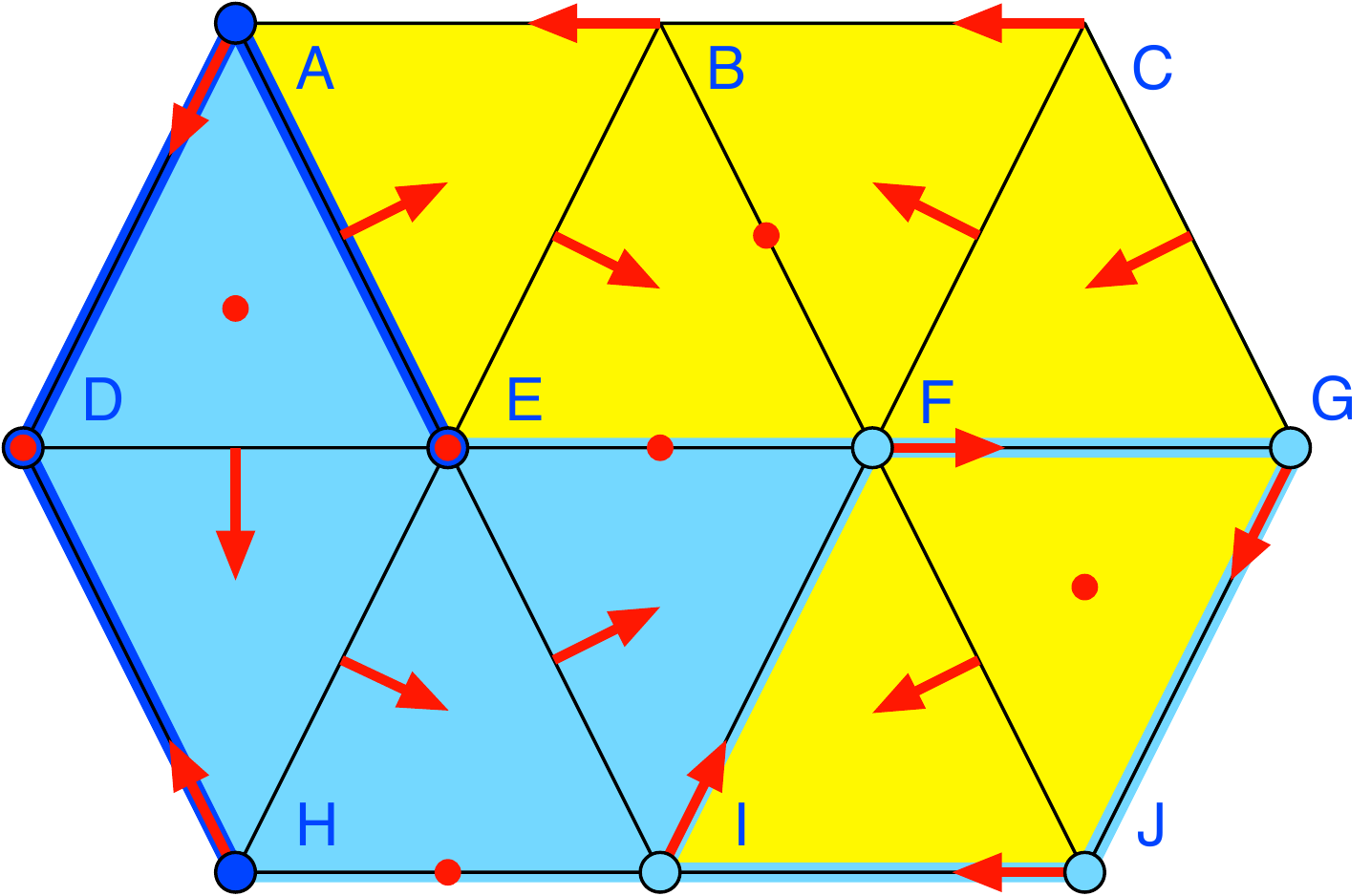}
  \caption{Sample isolated invariant set for the discrete vector
           field shown in Figure~\ref{fig:mainexample}. The
           simplices which belong to the isolated invariant
           set~$\cS$ are indicated in light blue, and are given
           by four vertices, nine edges, and four triangles.
           Its exit set~$\Exit\cS$ is shown in dark blue, and
           it consists of four vertices and three edges.}
  \label{fig:isolatedinvset}
\end{figure}

The two-dimensional example depicted in Figure~\ref{fig:mainexample} presents a simplicial complex which is built from~$10$
triangles, $19$ edges and~$10$ vertices, and a combinatorial vector field consisting of~$7$ critical cells
and a total of~$16$ vectors.
The set $\{ADE,DE,DEH,EF,EFI,EH,EHI,EI,F,FG,FI,G,GJ,HI,I,IJ,J\}$
is an example of an isolated invariant set for this combinatorial vector field.
It is presented in Figure~\ref{fig:isolatedinvset}.
Its exit set is $\{A,AD,AE,D,DH,E,H\}$ and its Conley index is the (co)homology of the pointed circle.
This isolated invariant set is not minimal. For instance, the singleton $\{EF\}$ is a subset which itself
is an isolated invariant set.

A {\em connection} from an isolated invariant set $\cS_1$ to an isolated invariant set $\cS_2$
is a sequence of vertices on a walk in $G_\cV$ originating in $\cS_1$ and terminating in $\cS_2$.
A family $\cM=\{\cM_p\,|\, p\in\PP\}$ indexed by a poset $\PP$
and consisting of mutually disjoint isolated invariant subsets of an isolated invariant set~$\cS$
is a {\em Morse decomposition} of~$\cS$ if any connection between elements in~$\cM$
which is not contained entirely
in one of the elements of $\cM$ originates in $\cM_{q'}$ and terminates in $\cM_{q}$ with $q'>q$.
The associated {\em Conley-Morse graph} is the partial order
induced on $\cM$ by the existence of connections, and represented as a directed graph labelled with
the Conley indices of the isolated invariant sets in $\cM$. Typically, the labels are written
as Poincar\'e polynomials, that is, polynomials whose $i$th coefficient equals the $i$th Betti
number of the Conley index.
\begin{figure}[tb]
  \includegraphics[width=0.99\textwidth]{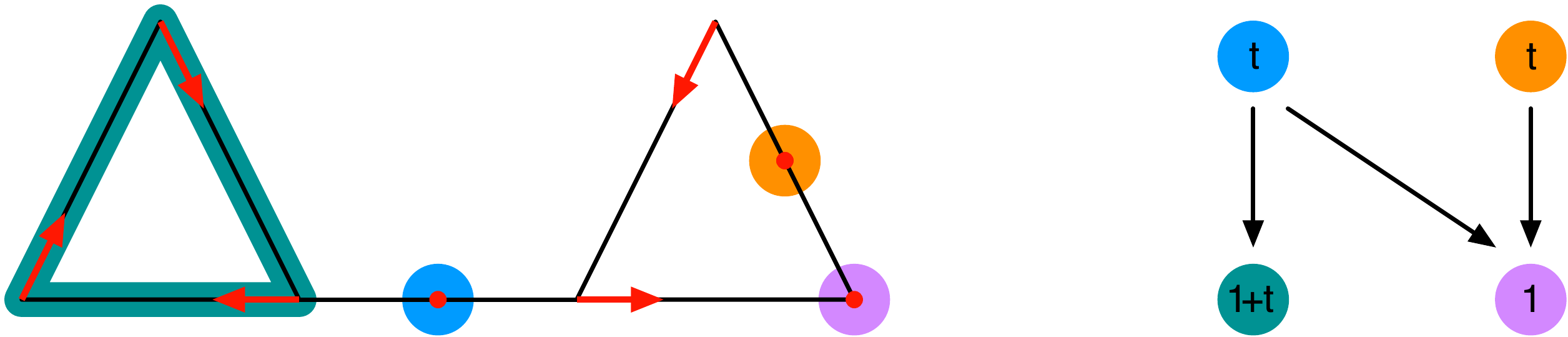}
  \caption{Morse decomposition for the example shown in
           Figure~\ref{fig:graphexample}. For this example,
           one can find four minimal Morse sets, which are
           indicated in the left image in different colors.
           The right image shows the associated Morse graph.}
  \label{fig:graphmorse}
\end{figure}

\begin{figure}[tb]
  \includegraphics[width=0.99\textwidth]{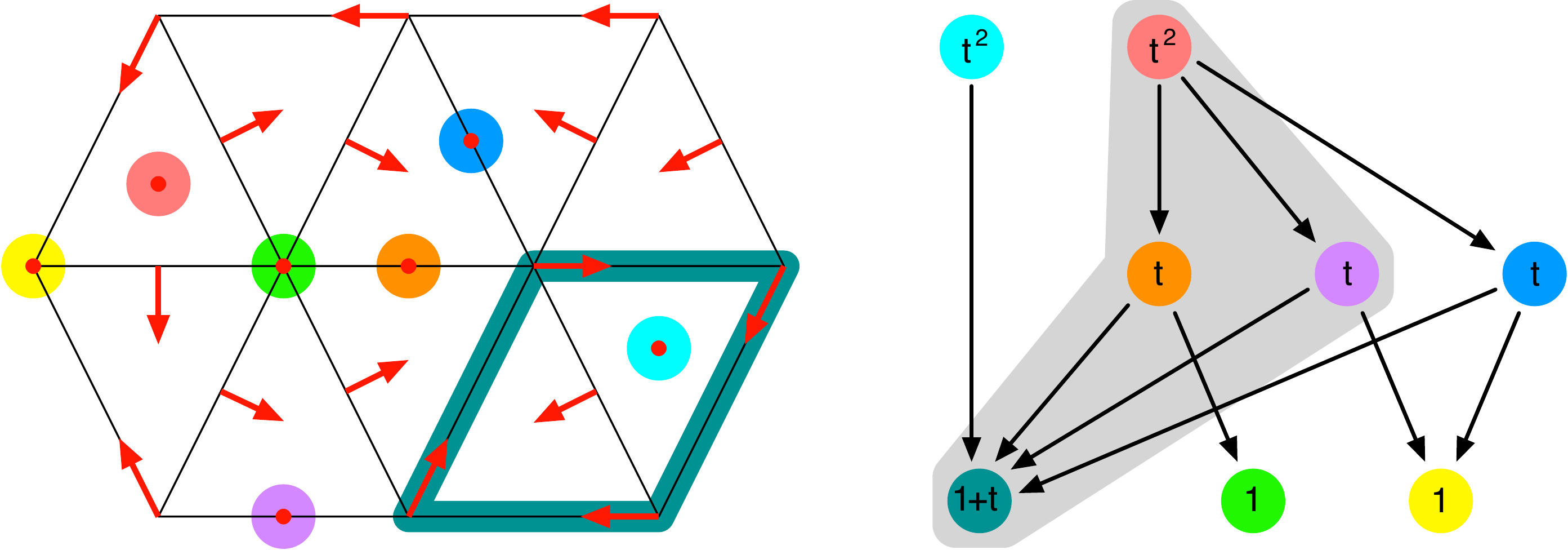}
  \caption{Morse decomposition for the example shown in
           Figure~\ref{fig:mainexample}. For this example,
           one can find eight minimal Morse sets, which are
           indicated in the left image in different colors.
           The right image shows the associated Morse graph.
           The isolated invariant set shown in
           Figure~\ref{fig:isolatedinvset} corresponds to
           the subgraph indicated by the gray shaded area
           in the Morse graph. }
  \label{fig:morsedecomp}
\end{figure}

An example of a Morse decomposition for the combinatorial vector field in Figure~\ref{fig:graphexample} is
\[
  \cM:=\{\{BF\}, \{F\}, \{DE\}, \{A,AD,C,CA,CD,D\}\} \; ,
\]
and the corresponding Conley-Morse graph is presented in Figure~\ref{fig:graphmorse}.
A Morse decomposition of the example in Figure~\ref{fig:mainexample} together with the associated
Conley-Morse graph is presented in Figure~\ref{fig:morsedecomp}.

The main result of this paper is the following theorem.

\begin{figure}[tb]
  \includegraphics[width=0.80\textwidth]{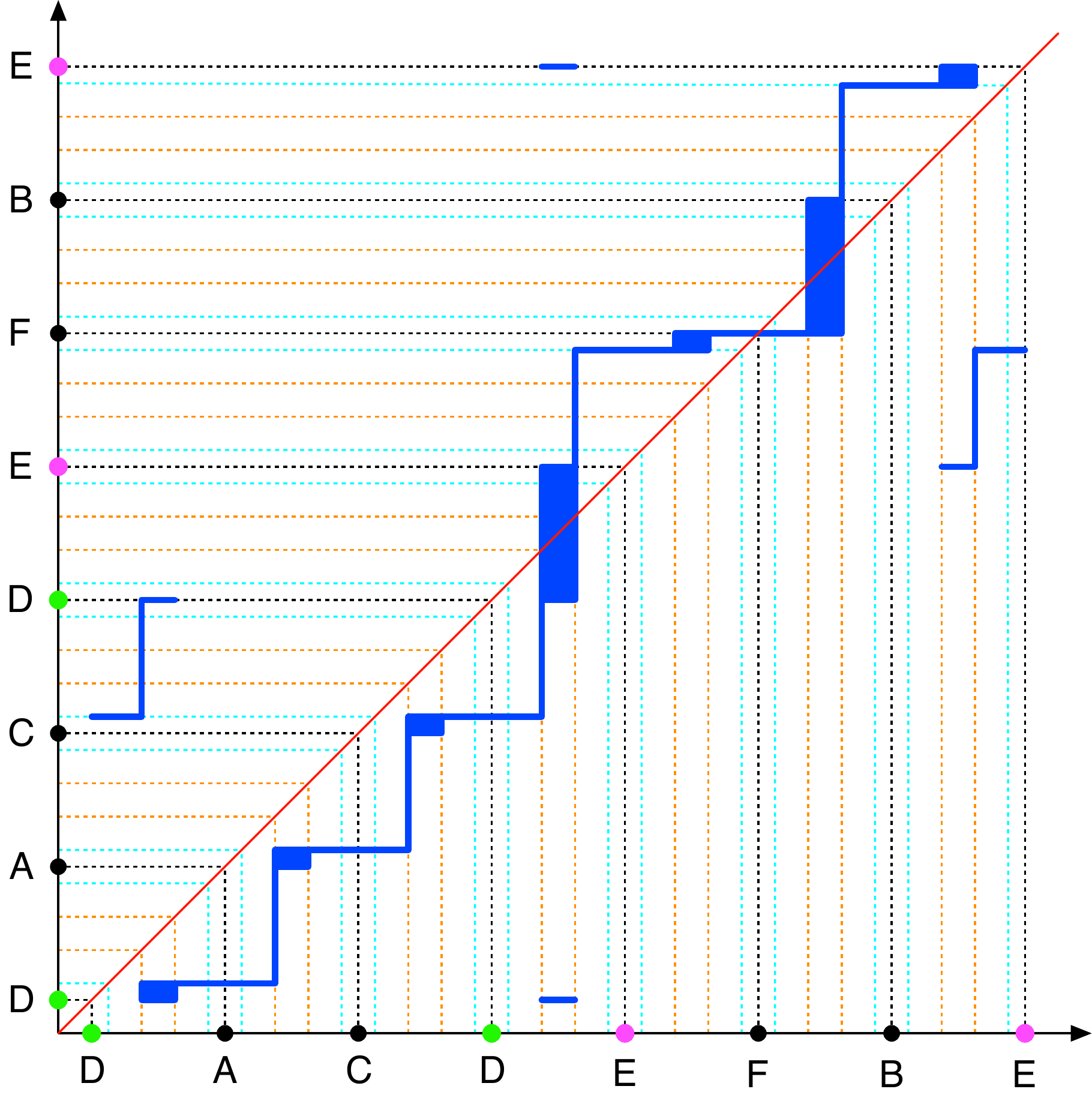}
  \caption{The multivalued map~$F$ for the combinatorial vector
           field shown in Figure~\ref{fig:graphexample}.
           For visualization purposes the domain of $F$ is straightened to
           a segment in which vertices $D$ (marked in green) and
           $E$ (marked in magenta) are represented twice.
           The graph of~$F$ is shown in blue. The edge $DE$ in the middle corresponds to the center
           edge in Figure~\ref{fig:graphexample}. To its left, the
           three line segments correspond to the cycle in the
           combinatorial vector field. Note that the two green vertices
           are identified.  The three edges to the right of the
           center correspond to the right triangle in
           Figure~\ref{fig:graphexample}. Also here the two magenta
           vertices are identified.}
  \label{fig:graphf}
\end{figure}

\begin{thm}
\label{thm:main}
For every combinatorial vector field $\cV$ on a simplicial complex $\cX$ there exists
an upper semicontinuous, acyclic, and inducing identity in homology multivalued map $F:|\cX|\mto|\cX|$
on the geometric realization $|\cX|$ of $\cX$ such that
\begin{itemize}
   \item[(i)] for every Morse decomposition~$\cM$
   of~$\cV$ there exists a Morse decomposition~$M$ of the semidynamical system induced by~$F$,
   \item[(ii)] the Conley-Morse graph of~$M$ is
   isomorphic to the Conley-Morse graph of~$\cM$,
   \item[(iii)] each element of~$M$ is contained in the geometric representation of the corresponding element of~$\cM$.
\end{itemize}
\end{thm}
This theorem is an immediate consequence of the much more detailed theorems presented
in Section~\ref{sec:correspondence}.
The multivalued map $F$ guaranteed by Theorem~\ref{thm:main} for the example in Figure~\ref{fig:graphexample}
is presented in Figure~\ref{fig:graphf}.

\section{Conley theory for multivalued topological dynamics}
\label{sec:conley}

In this section we recall the main concepts of Conley theory for multivalued dynamics
in the combinatorial and classical setting:
isolated invariant sets, index pairs, Conley index and Morse decompositions.

\subsection{Preliminaries}
\label{sec:pre}
We write $f: X\omap Y$ to denote a {\em partial function}, that is, a function
whose {\em domain}, denoted $\dom f$, is a subset of $X$.
We write $\im f:=f(X)$ to denote the {\em image} of $f$
and $\Fix f:=\setof{x\in\dom f\mid f(x)=x}$ to denote the set of {\em fixed points} of $f$.

Given a topological space $X$ and a subset $A\subset X$,
we denote by $\cl A$, $\Int A$ and $\bd A$ respectively
the {\em closure}, the {\em interior} and the {\em boundary} of $A$.
We often use the set $\exit A:=\cl  A \setminus  A$ which we call the {\em exit set} or {\em mouth} of $ A$.
Whenever applying an operator like $\cl$ or $\exit$ to a singleton, we drop the braces to keep the notation simple.

The {\em singular cohomology} of the pair  $(X,A)$ is denoted $H^*(X,A)$.
Note that in this paper we apply cohomology only to polyhedral  pairs or pairs
weakly homotopy equivalent
to polyhedral pairs. Hence, the singular cohomology is the same as Alexander-Spanier
cohomology. In particular, all but a finite number of Betti numbers of the pair $(X,A)$ are zero.
The corresponding {\em Poincar\'e polynomial} is the polynomial whose $i$th coefficient is the $i$th Betti number.

By a multivalued map $F:X\mto X$ we mean a map from $X$ to the family of non-empty subsets of $X$.
We say that $F$ is {\em upper semicontinuous} if for any open $U\subset X$  the set
$\{x\in X\mid  F(x)\subset U\}$ is open.
We say that $F$ is {\em strongly upper semicontinuous} if for any $x\in X$ there
exists a neighborhood $U$ of $X$ such that $x'\in U$ implies $F(x')\subset F(x)$.
Note that every strongly upper semicontinuous multivalued map is upper semicontinuous.
We say that $F$ is {\em acyclic-valued} if $F(x)$ is acyclic for any $x\in X$.

We consider a {\em simplicial complex} as a finite family $\cX$ of finite sets such that
any non-empty subset of a set in $\cX$ is in $\cX$.
We refer to the elements of $\cX$ as simplices.
By the {\em dimension} of a simplex we mean one less than its cardinality.
We denote by $\cX_k$ the set of simplices of dimension $k$.
A {\em vertex} is a simplex of dimension zero.
If $\sigma,\tau\in\cX$ are simplices and $\tau\subset\sigma$ then we say that $\tau$ is a {\em face} of $\sigma$
and $\sigma$ is a {\em coface} of $\tau$.
An $(n-1)$-dimensional face of an $n$-dimensional simplex is called
a {\em facet}.
We say that a subset $\cA\subset \cX$ is {\em open} if
all cofaces of any element of $\cA$ are also in $\cA$.
It is easy to see that the family
of all open sets of $\cX$ is a $T_0$ topology on $\cX$, called {\em Alexandroff topology}.
It corresponds to the face poset of $\cX$ via the Alexandroff Theorem \cite{Al1937}.
In particular,  the {\em closure} of $\cA\subset \cX$ in the Alexandroff topology
consists of all faces of simplices in $\cA$.
To avoid confusion, in the case of Alexandroff topology we write $\ccl \cA$ and $\Exit \cA$
for the closure and the exit of $\cA\subset\cX$.

By identifying  vertices of an $n$-dimensional simplex $\sigma$ with a collection of
 $n+1$ linearly independent vectors in $\RR^d$ with $d>n$ we obtain a {\em geometric realization} of $\sigma$.
 We denote it by $|\sigma|$. However,  whenever the meaning is clear from the context
 we drop the bars to keep the notation simple.
By choosing the identification in such a way that all vectors corresponding to vertices of $\cX$
are linearly independent we obtain a {\em geometric realization} of $\cX$ given by
\[
   |\cX|:=\bigcup_{\sigma\in\cX}|\sigma|.
\]
Note that up to a homeomorphism the geometric realization
does not depend on a particular choice of the identification.
In the sequel we assume that a simplicial complex $\cX$ and its geometric realization $X:=|\cX|$ are fixed.
Given a vertex $v\in\cX_0$, we denote by $t_v:|\cX|\to[0,1]$ the map
which assigns to each point~$x\in|\cX|$ its barycentric coordinate
with respect to the vertex~$v$. For a simplex $\sigma\in\cX$ the {\em open cell} of $\sigma$ is
\[
    \cell\sigma:=\setof{x\in|\sigma|\mid t_v(x)>0 \text{ for } v\in\sigma}.
\]
For $\cA\subset\cX$ we write
\[
  \ogr{\cA}:=\bigcup_{\sigma\in\cA}\cell\sigma.
\]
One easily verifies the following proposition.
\begin{prop}
\label{prop:ogr}
We have the following properties
\begin{itemize}
   \item[(i)] if $\cA$ is closed in $\cX$ then $|\cA|=\ogr{\cA}$,
   \item[(ii)] if $\Exit\cA$ is closed in $\cX$ then $\ogr{\cA}=|\cA|\setminus|\Exit\cA|$.
\end{itemize}
\qed
\end{prop}


\subsection{Combinatorial case}
\label{sec:comb-dyn}

The concept of a combinatorial vector field was introduced by Forman~\cite{Fo98}.
There are a few equivalent ways of stating its definition.
The definition introduced in Section \ref{sec:main-result} is among the simplest:
a {\em combinatorial vector field} on a simplicial complex $\cX$ is a partition $\cV$ of $\cX$
into singletons and doubletons such that each doubleton consists of a simplex and
one of its facets. The partition induces
an injective partial map which sends the element of each singleton to itself and each facet
in a doubleton to its coface in the same doubleton.
This leads to the following equivalent definition which will be used in the rest of the paper.

\begin{defn} \label{def:discrete:vector:field} (see \cite[Definition 3.1]{KaMrWa16})
An injective partial
self-map $\cV:\cX \nrightarrow \cX$ of a simplicial complex $\cX$ is called a {\em combinatorial vector
field}, or also a {\em discrete vector field} if
\begin{itemize}
  \item[(i)]   For every simplex $\sigma \in \dom\cV$ either
               $\cV(\sigma) = \sigma$, or~$\sigma$ is a facet of~$\cV(\sigma)$.
  \item[(ii)]  $\dom\cV \cup \im\cV = \cX$,
  \item[(iii)]  $\dom\cV \cap \im\cV = \Fix \cV$.
\end{itemize}
\end{defn}
Note that every combinatorial vector field is a special case of a combinatorial multivector field
introduced and studied in \cite{Mr15}.

Given   a combinatorial
vector field $\cV$ on $\cX$, we define the associated {\em combinatorial multivalued
flow} as the multivalued map $\Pi_\cV : \cX \mto \cX$ given by
\begin{equation} \label{def:discrete:flow}
   \Pi_\cV(\sigma) :=
   \begin{cases}
     \ccl\sigma & \text{if $\sigma\in\Fix\cV$} \; , \\
     \Bd\sigma \setminus\{ \cV^{-1}(\sigma) \}
     & \text{if $\sigma \in \im\cV \setminus \Fix(\cV)$} \; ,\\
     \{\cV(\sigma)\} & \text{if $\sigma \in \dom\cV \setminus
       \Fix(\cV)$} \; .
   \end{cases}
\end{equation}
For the rest of the paper we assume that
$\cV$ is a fixed combinatorial vector field on $\cX$ and $\Pi_\cV$ denotes
the associated combinatorial multivalued flow.

A {\em solution} of the flow~$\Pi_\cV$ is a partial
function $\rho : \ZZ \nrightarrow \cX$ such that
$\rho(i+1) \in \Pi_\cV(\rho(i))$ whenever  $i, i+1 \in \dom\rho$.
The solution $\rho$ is {\em full} if $\dom\rho=\ZZ$.
The {\em invariant part} of $\cS \subset \cX$, denoted $\Inv\cS$, is the collection
of those simplices $\sigma \in \cS$ for which there exists a full solution~$\rho : \ZZ \to \cS$
such that $\rho(0) = \sigma$.
A set $\cS \subset \cX$ is  {\em invariant\/} if $\Inv\cS=\cS$.

\begin{defn} \label{def:discrete-isolated-inv} (see \cite[Definition 3.4]{KaMrWa16})
A subset $\cS \subset \cX$, invariant with respect to a combinatorial vector field $\cV$,
is called an {\em isolated invariant set\/} if the {\em exit set}
$
  \Exit\cS = \ccl \cS \setminus \cS
$
is closed and there is no solution $\rho : \{-1,0,1\}
\to \cX$ such that $\rho(-1), \rho(1) \in \cS$  and $\rho(0) \in \Exit\cS$.
The closure~$\ccl\cS$ is called
an {\em isolating block\/} for the isolated invariant
set~$\cS$.
\end{defn}

\begin{prop} \label{prop:exSclosed+-} (see \cite[Proposition 3.7]{KaMrWa16})
An invariant set~$\cS\subset \cX$ is an isolated invariant set if ~$\Exit~\cS$ is closed
and for every $\sigma \in \cX$ we have
$\sigma^- \in \cS$ if and only if $\sigma^+ \in \cS$, where
\begin{equation} \label{def:sigma:pm}
  \sigma^+ := \begin{cases}
                \cV(\sigma) & \text{if $\sigma\in\dom\cV$} \\
                \sigma & \text{otherwise}
              \end{cases}
  \qquad\mbox{and}\qquad
  \sigma^- := \begin{cases}
                \sigma & \text{if $\sigma\in\dom\cV$}\\
                \cV^{-1}(\sigma) & \text{otherwise}
              \end{cases} .
\end{equation}
\end{prop}
An immediate consequence of Proposition~\ref{prop:exSclosed+-} is the following corollary.
\begin{cor}
\label{cor:convexity}
If $\cS$ is an isolated invariant set then for any $\tau\in\cX$ and $\sigma,\sigma'\in\cS$ we have
\[
     \sigma\subset\tau\subset\sigma'  \implies \tau\in\cS.
\]\qed
\end{cor}

A pair $\cP=(\cP_1, \cP_2)$ of closed subsets of $\cX$ such that $\cP_2\subset \cP_1$
is an {\em index pair} for $\cS$ if the following three conditions are satisfied
\begin{gather}
  \cP_1\cap\Pi_{\cV}(\cP_2)\subset \cP_2,\label{eq:ip1}\\
  \Pi_{\cV}(\cP_1\setminus\cP_2)\subset \cP_1,\label{eq:ip2}\\
  \cS=\Inv( \cP_1\setminus \cP_2).\label{eq:ip3}
\end{gather}

By \cite[Theorem 7.11]{Mr15} the pair $(\cl \cS,\Exit \cS)$ is an index pair for
$\cS$ and the (co)homology of the index pair of $\cS$ does not depend on the particular choice of index pair
but only on $\cS$. Hence, by definition, it is the {\em Conley index} of $\cS$.
We denote it $\Con(\cS)$.


\subsection{Classical case}
\label{sec:class-dyn}

The study of the Conley index for multivalued maps was initiated in \cite{KaMr95}
with a restrictive concept of the isolating neighborhood,  limiting possible applications.
In particular, that theory is not satisfactory for the needs of this paper.
These limitations were removed by a new theory developed recently in \cite{BaMr16,Ba17}.
We recall the main concepts of the generalized theory below.

Let $F:X\mto X$ be an upper semicontinuous map with compact, acyclic values.
A partial map $\rho:\ZZ\omap X$ is called a {\em solution for $F$ through $x\in X$} if we have both $\rho (0)=x$ and $\rho (n+1)\in F(\rho(n))$ for all $n,n+1\in \dom\rho$.
Given $N\subset X$ we define its {\em invariant part} by
\[
\Inv N:=\{x\in N\mid  \exists\,\rho:\ZZ\to N\mbox{ which is a solution  for }F\mbox{ through }x\}.
\]
A compact set $N\subset X$ is an {\em isolating neighborhood} for $F$ if $\Inv N \subset \inte N$.
The $F$-{\em boundary} of a given set $A\subset X$ is
$$
\bd _F(A):=\cl A\cap\cl (F(A)\setminus A).
$$

\begin{defn}
{\em
\label{def:weak-index-pair}
A pair $P=(P_1,P_2)$ of compact sets $P_2\subset P_1 \subset N$ is a {\em weak index pair} for $F$ in $N$ if the following properties are satisfied.
\begin{aenum}
\item $F(P_i)\cap N \subset P_i$ for $i=1,2$\;,
\item $\bd_F(P_1)\subset P_2$,
\item $\Inv N \subset \Int(P_1\setminus P_2)$,
\item $P_1\setminus P_2\subset \Int N$.
\end{aenum}
}
\end{defn}

For the weak index pair $P$ we set
\[
T(P):=T_N(P):=(P_1\cup (X\setminus \Int N), P_2\cup (X\setminus \Int N)).
\]
and define the associated index map $I_P$ as the composition $H^*(F_P)\circ H^*(i_P)^{-1}$, where
$F_P:P\mto T(P)$ is the restriction of $F$ and $i_P:P\mto T(P)$ is the inclusion map.
The module $\ClasCon(S,F):=L(H^*(P),I_P)$, where $L$ is the Leray functor (see \cite{Mr90a})
is the {\em cohomological Conley index} of the isolated invariant set $S$.
The correctness of the definition is the consequence of the following two results.

\begin{thm} (see \cite[Theorem 4.12]{BaMr16})
\label{thm:existence}
For every neighborhood $W$ of $\Inv N$ there exists a weak index pair $P$
in $N$ such that $P_1\setminus P_2\subset W$.
\end{thm}

\begin{thm}\label{th_ind} (see \cite[Theorem 6.4]{BaMr16})
The  module $L(H^*(P),I_P)$ is independent of the choice of an isolating neighborhood $N$ for $S$ and of a weak index pair~$P$ in~$N$.
\end{thm}

\subsection{Morse decompositions.}

In order to formulate the definition of the Morse decomposition of an isolated invariant set
we need the concepts of $\alpha$- and $\omega$-limit sets. We formulate theses definitions independently
in the combinatorial and classical settings.
Given a full solution $\rho:\ZZ\to\cX$ of the combinatorial dynamics $\Pi_\cV$ on $\cX$,
the {\em $\alpha$- and $\omega$-limit sets} of $\rho$
are respectively the sets
\[
   \alpha(\rho):=\bigcap_{n\in\ZZ}\{\rho(k)\mid k\geq n\}, \quad\quad
   \omega(\rho):=\bigcap_{n\in\ZZ}\{\rho(k)\mid k\leq n\}.
\]
Note that $\alpha$- and $\omega$-limit sets of $\Pi_\cV$ are always
non-empty invariant sets, because~$\cX$ is finite.

Now, given a solution $\phi:\ZZ\to X$ of a multivalued upper semicontinuous map $F:X \mto X$,
we define its {\em $\alpha$- and $\omega$-limit sets} respectively by
$$
\alpha(\phi):=\bigcap_{k\in\ZZ}\cl\phi((-\infty,-k])  \quad\quad  \omega(\phi):=\bigcap_{k\in\ZZ}\cl\phi([k,+\infty)).
$$

\begin{defn}\label{def:Morse_dec}
{\rm
Let $\cS$ be an isolated invariant set of $\Pi_\cV:\cX\mto \cX$.
We say that the family $\cM:=\{\cM_r\,|\,r\in\PP\}$ indexed by a poset $\PP$
is a {\em Morse decomposition of $\cS$} if  the following conditions are satisfied:
\begin{aenum}
\item the elements of $\cM$ are mutually disjoint isolated invariant subsets of $\cS$,
\item for every full solution $\phi$ in $X$ there exist $r,r'\in\PP$, $r\leq r'$, such that $\alpha(\phi)\subset \cM_{r'}$ and $\omega(\phi)\subset \cM_r$,
\item if for a full solution $\phi$ in $\cX$ and $r\in\PP$ we have $\alpha(\phi)\cup\omega(\phi)\subset \cM_r$, then $\im\phi\subset \cM_r$.
\end{aenum}
By replacing in the above definition the multivalued map $\cF$ on $\cX$ by an upper semicontinuous map $F:X\to X$
and adjusting the notation accordingly we obtain the definition of the {\em Morse decomposition} $M:=\{M_r\,|\,r\in\PP\}$
of an isolated invariant set $S$ of $F$.
}
\end{defn}
It is not difficult to observe that  Definition~\ref{def:Morse_dec} in the combinatorial setting
is equivalent to the brief definition of Morse decomposition given in terms of connections
in Section~\ref{sec:main-result}. Moreover, in the case of combinatorial vector fields
Definition~\ref{def:Morse_dec} coincides with the definition presented in \cite[Section 9.1]{Mr15}.


\section{From combinatorial to classical dynamics}
\label{sec:map-F}

In this section, given a combinatorial vector field $\cV$ on a simplicial complex $\cX$,
we recall from \cite{KaMrWa16} the construction of a multivalued self-map $F=F_\cV:X\mto X$
on the geometric realization $X:=|\cX|$ of $\cX$. This map will be used
to establish the correspondence of Conley indices, Morse decompositions and Conley-Morse graphs
between the combinatorial and classical multivalued dynamics.

\subsection{Cellular decomposition}
We begin by recalling a special
cellular complex representation of~$X = |\cX|$ used in the
construction of the multivalued map~$F$. For this we need some terminology.
Let $d$ denote the maximal dimension of the simplices in $\cX$.
Fix a $\lambda\in\RR$ such that $0 \le \lambda < \frac{1}{d+1}$ and a point $x\in X$.
The {\em $\lambda$-signature} of~$x$ is the function
\begin{equation} \label{lambda:signature}
   \sign^\lambda x : \cX_0\ni v\mapsto \sgn\left( t_v(x) - \lambda \right)\in \{-1,0,1\},
\end{equation}
where $\sgn : \RR \to \{-1,0,1\}$ is the standard sign function.
Then a simplex $\sigma \in \cX$ is a
{\em $\lambda$-characteristic simplex\/} of~$x$ if both
$\sign^\lambda x|_{\sigma} \geq 0$
and
$(\sign^\lambda x)^{-1}(\{1\}) \subset \sigma$ are satisfied.
We denote the family of $\lambda$-characteristic simplices of~$x$ by
\begin{displaymath}
  \cX^\lambda(x) :=
  \setof{\sigma\in \cX \; \mid \; (\sign^\lambda x)^{-1}(\{1\})\subset \sigma
    \;\text{ and }\; \sign^\lambda x(v)\geq 0 \;\text{ for all }\;
    v\in \sigma}
  \; .
\end{displaymath}
For any $\lambda \geq 0$, the set $(\sign^\lambda x)^{-1}(\{1\})$ is a
simplex. We call it the {\em minimal characteristic
simplex} of~$x$ and we denote it
by~$\sigma^\lambda_{min}(x)$.
Note that
\begin{equation}
\label{eq:min-cell-sigma}
\sigma=\sigma^0_{min}(x) \iff x\in\cell\sigma.
\end{equation}
If $\lambda > 0$, then the
set $(\sign^\lambda x)^{-1}(\{0,1\})$ is also a simplex.
We call it
the {\em maximal characteristic simplex} of~$x$ and we
denote it by~$\sigma^\lambda_{max}(x)$.

\begin{lem} \label{lem:epsilon-delta-include}  (see \cite[Lemma 4.2]{KaMrWa16})
If $0 \le \epsilon < \lambda < \frac{1}{1+d}$, then $\sigma^\lambda_{max}(x)
\subset \sigma^\epsilon_{min}(x)$ for any $x \in X = |\cX|$.
\end{lem}

\begin{figure}[tb]
  \includegraphics[width=0.80\textwidth]{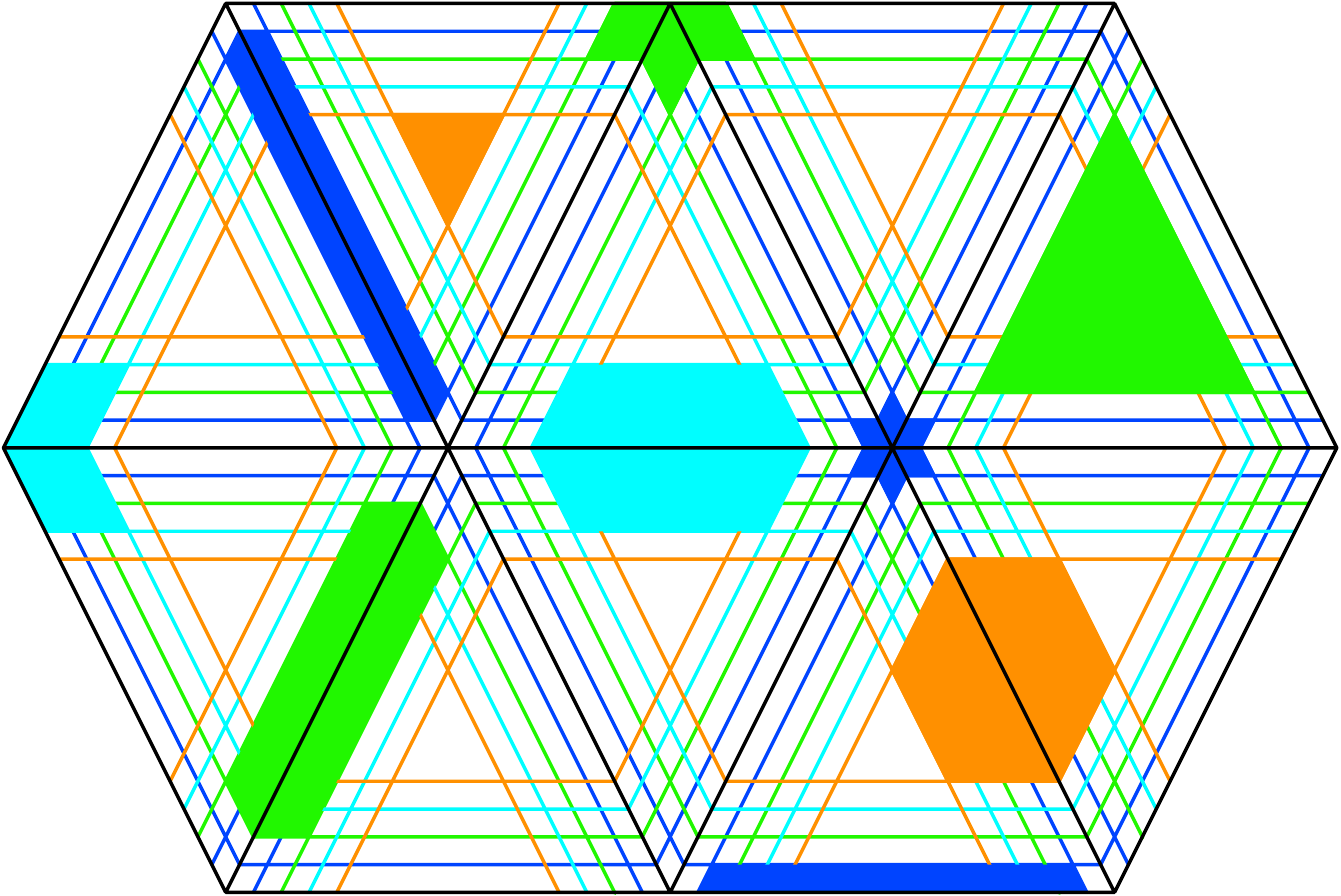}
  \caption{Sample cell decomposition boundaries for the simplicial
           complex~$\cX$ from Figure~\ref{fig:mainexample}. The
           colored lines indicate the boundaries of $\epsilon$-cells
           (orange), $\gamma$-cells (cyan), $\delta$-cells
           (green), and $\delta'$-cells (blue). Throughout the
           paper, we assume that $0 < \delta' < \delta < \gamma
           < \epsilon$. The figure also contains ten sample cells:
           Two orange $\epsilon$-cells which are associated with
           a $2$-simplex (upper left) and a $1$-simplex (lower
           right); two cyan $\gamma$-cells which correspond to
           a $1$-simplex (middle) and a $0$-simplex (left);
           three green $\delta$-cells for a $2$-simplex (upper right),
           a $1$-simplex (lower left), and a $0$-simplex (top
           middle); as well as three blue $\delta'$-cells for two
           $1$-simplices (upper left and bottom right) and a
           $0$-simplex (right middle). All of these cells are
           open subsets of~$|\cX|$.}
  \label{fig:celldecomp}
\end{figure}

\begin{lem} \label{lem:cs-semicontiniuty} (see \cite[Lemma 4.3]{KaMrWa16})
For all $x \in X = |\cX|$ we have~$\cX^\lambda(x)\neq\emptyset$.
Moreover, there exists a
neighborhood~$U$ of the point~$x$ such that
$\cX^\lambda(y) \subset \cX^\lambda(x)$   for all $ y \in U$.
\end{lem}

Given a $\sigma \in \cX$  by a
{\em $\lambda$-cell generated by~$\sigma$} we mean
\begin{displaymath}
  \csl{\sigma} :=
  \left\{ x \in X \; \mid \;
    \cX^\lambda(x) = \{ \sigma \} \right\}.
\end{displaymath}
We recall (cf. \cite[Formulas~(12) and~(13)]{KaMrWa16})  the following characterizations
of  $\csl{\sigma}$
and its closure in terms of barycentric coordinates:
\begin{eqnarray}\label{eq:sigma_l}
\csl{\sigma} &=& \left\{x\in X \; \mid \;
    t_v(x) > \lambda \;\mbox{ for  }\; v \in \sigma
    \;\mbox{ and }\;
    t_v(x) < \lambda
    \;\mbox{ for  }\; v \notin \sigma \right\},\\
\label{eq:cl_sigma_l}
\cl\csl{\sigma} &=& \left\{x\in X \; \mid \;
    t_v(x) \ge \lambda \;\mbox{ for  }\; v \in \sigma
    \;\mbox{ and }\;
    t_v(x) \le \lambda
    \;\mbox{ for  }\; v \notin \sigma \right\}.
\end{eqnarray}
Then the following proposition follows easily from \eqref{eq:sigma_l}.
\begin{prop}
\label{prop:sigma-decomp}
For $\lambda$ satisfying $0<\lambda<\frac{1}{d+1}$
the $\lambda$-cells are open in $|\cX|$ and mutually disjoint.
\qed
\end{prop}

Another characterization of $\cl\csl{\sigma}$ is given by the following corollary.
\begin{cor} \label{cor:equiv-cell-closure} (see \cite[Corollary 4.6]{KaMrWa16})
The following statements are equivalent:
\begin{itemize}
\item[(i)] $\sigma \in\cX^\lambda(x)$,
\item[(ii)] $\sigma^\lambda_{min}(x) \subset \sigma
\subset \sigma^\lambda_{max}(x)$,
\item[(iii)] $x \in \cl\csl{\sigma}$.
\end{itemize}
\end{cor}
%

Note that for a simplicial complex $\cX$ we have the following easy to verify formula:
\begin{equation}
\label{eq:x}
|\cX|=\bigcup_{\sigma\in\cX}\cl\csd{\sigma}.
\end{equation}

The cells $\csl{\sigma}$ for various values of $\lambda$ are visualized in Figure~\ref{fig:celldecomp}.
They are the building blocks for the multivalued map $F$.

\subsection{The maps $F_\sigma$ and the map $F$}
\label{sec:Fsigma-F}

We now recall from \cite{KaMrWa16} the construction of the strongly upper semi-continuous
map~$F$ associated with a combinatorial vector field.
For this, we fix two constants
\begin{equation} \label{def:gammaepsbound}
  0 < \gamma < \epsilon < \frac{1}{d+1}
\end{equation}
and for any $\sigma \in \cX$ we set
\begin{eqnarray}
  A_\sigma & := & \left\{ x \in \sigma^+ \; \mid \;
    t_v(x) \geq \gamma \;\mbox{ for all }\; v \in \sigma^- \right\}
    \; \cup \; \sigma^- \; , \nonumber \\[0.5ex]
  B_\sigma & := & \left\{ x \in \sigma^+ \; \mid \;
    \mbox{there exists a }\; v \in \sigma^- \;\mbox{ with }\;
    t_v(x) \leq \gamma \right\} \; ,
    \label{def:ABCsigma} \\[0.5ex]
  C_\sigma & := & A_\sigma \cap B_\sigma \; . \nonumber
\end{eqnarray}
Then the following lemma is an immediate consequence of \cite[Lemma 4.8]{KaMrWa16}.
\begin{lem} \label{lem:ABC-contractible}
For any simplex
$\sigma \in \cX \setminus \Fix\cV$ the sets $A_\sigma$, $B_\sigma$ and $C_\sigma$
are contractible.
\end{lem}

\medskip

For every simplex~$\sigma \in \cX$ we define a multivalued map~$F_\sigma : X
\mto X$ by

\begin{equation} \label{def:Fsigma}
  F_\sigma (x) :=
  \begin{cases}
    \emptyset & \mbox{if }\; \sigma\notin \cX^\epsilon(x)
      \; , \\[0.5ex]
    A_\sigma  & \mbox{if }\; \sigma \in \cX^\epsilon (x)
      \; , \;\; \sigma \neq \sigma^\epsilon_{max}(x)^+ \; ,
      \; \mbox{ and } \; \sigma \neq \sigma^\epsilon_{max}(x)^-
      \; , \\[0.5ex]
    B_\sigma & \mbox{if }\; \sigma = \sigma^\epsilon_{max}(x)^+
      \neq \sigma^\epsilon_{max}(x)^- \; , \\[0.5ex]
    C_\sigma & \mbox{if }\; \sigma = \sigma^\epsilon_{max}(x)^-
      \neq \sigma^\epsilon_{max}(x)^+ \; , \\[0.5ex]
    \sigma   & \mbox{if }\; \sigma = \sigma^\epsilon_{max}(x)^-
      = \sigma^\epsilon_{max}(x)^+ \; ,
  \end{cases}
\end{equation}
and the multivalued map
$F : X \mto X$ by

\begin{equation}\label{eq:defF}
  F(x) := \bigcup_{\sigma \in \cX} F_\sigma(x)
  \qquad\mbox{ for all }\qquad
  x \in X = |\cX| \; .
\end{equation}
Figure~\ref{fig:graphf} shows the graph $\{(x,y)\in X\times X\mid y\in F(x)\}$ of
the so-constructed map~$F$ for the vector field in Figure~\ref{fig:graphexample}.

\medskip

One of main results proved in \cite{KaMrWa16} is the following theorem.

\begin{thm} \label{th:F-usc-acyclic} (see \cite[Theorem 4.12]{KaMrWa16})
The map~$F$ is strongly upper semicontinuous
and for every $x\in X$ the set~$F(x)$ is non-empty
and contractible.
\qed
\end{thm}

%
%
%
\section{The correspondence between combinatorial and classical dynamics}
\label{sec:correspondence}

In this section we present the constructions and theorems establishing the correspondence between the multivalued dynamics of
a combinatorial vector field $\cV$ on the simplicial complex $\cX$ and the associated multivalued dynamics of
the multivalued map $F=F_\cV$ constructed in Section~\ref{sec:map-F}. The theorems presented in this section provide
the proof of Theorem~\ref{thm:main}.

Throughout the section we assume that $d$ is the maximal
dimension of the simplices in $\cX$.

\subsection{Correspondence of isolated invariant sets.}
\begin{figure}[tb]
  \includegraphics[width=0.60\textwidth]{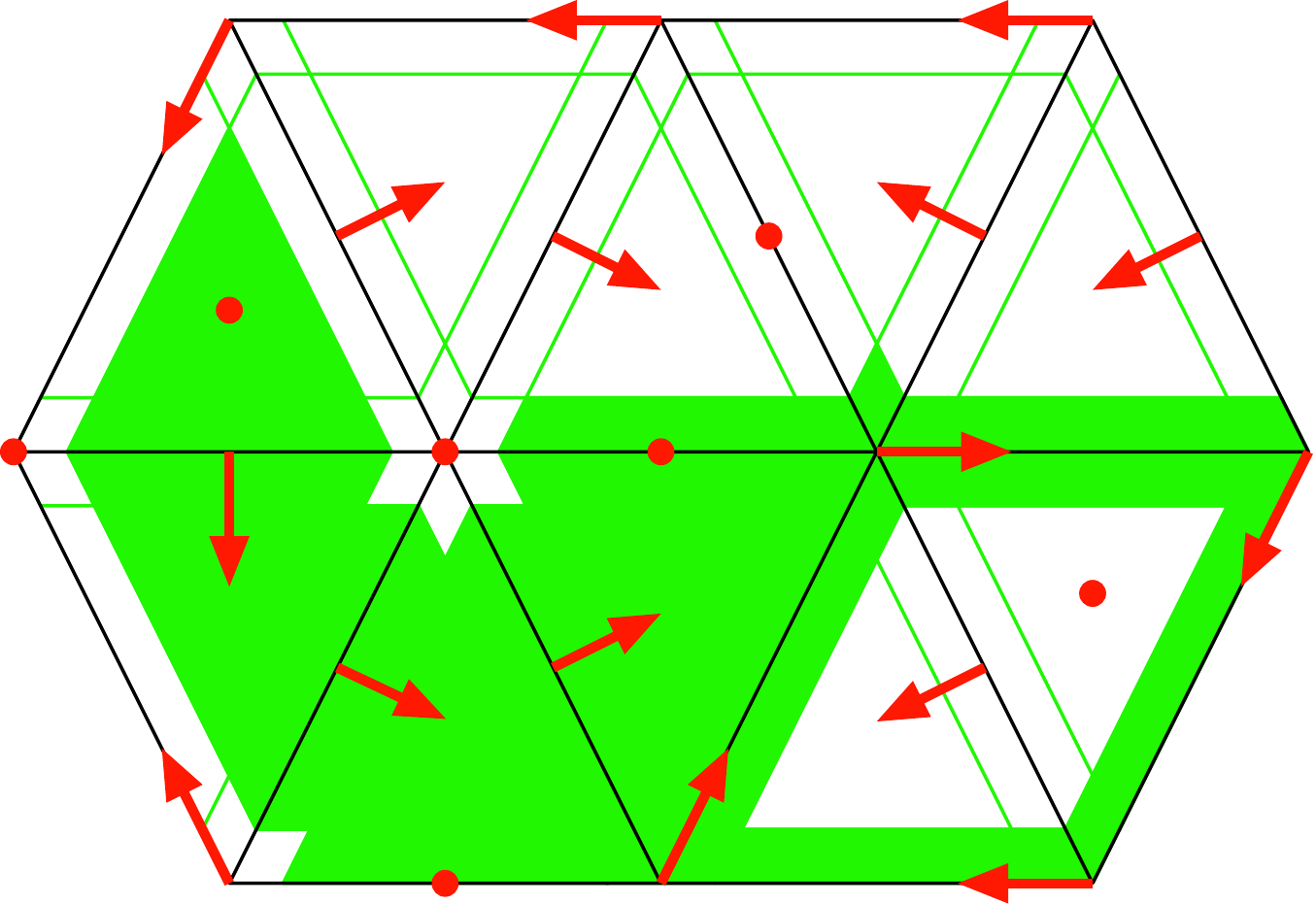}
  \caption{Isolating block~$N_\delta$ for the isolated invariant
           set~$\cS$ shown in Figure~\ref{fig:isolatedinvset}.
           Notice that the block is the union of closed
           $\delta$-cells. For reference, we also show the
           $\delta$-cell boundaries outside~$N_\delta$, but
           these are not part of the isolating block.
           The block is homeomorphic to a closed annulus.}
  \label{fig:isoblock}
\end{figure}

In order to establish the correspondence on the level of isolated invariant sets we fix a constant $\delta$
satisfying
\[
  0  < \delta < \gamma < \epsilon < \frac{1}{d+1} \; ,
\]
where $\gamma$ and $\epsilon$ are the constants chosen in Section~\ref{sec:Fsigma-F} (see \eqref{def:gammaepsbound}).
For $\cA\subset\cX$ and any constant~$\beta$ satisfying $0<\beta<\frac{1}{d+1}$
we further set
\begin{equation}
\label{eq:N-beta}
N_\beta(\cA):= \bigcup_{\sigma \in \cA} \cl\csb{\sigma}.
\end{equation}
Let $\cS\subset\cX$ be an isolated invariant set for the combinatorial vector field $\cV$
in the sense of Definition~\ref{def:discrete-isolated-inv}.
The following theorem associates with~$\cS$ an isolating block for~$F$, and it was proved in \cite{KaMrWa16}.
\begin{thm} \label{th:N-iso-block} (see \cite[Theorem 5.7]{KaMrWa16})
The set
\begin{equation} \label{eq:N-delta}
 N :=N_\delta:= N_\delta(\cS)
\end{equation}
is an isolating block for~$F$. In particular, it is an isolating neighborhood for $F$.\qed
\end{thm}
A sample of an isolating block for the map $F$ given by (\ref{eq:defF})
which corresponds to the combinatorial isolated invariant set in Figure~\ref{fig:isolatedinvset}
is presented in Figure~\ref{fig:isoblock}.

Theorem~\ref{th:N-iso-block} lets us associate with $\cS$ an isolated invariant set
\[
\isoS:=\Inv N_\delta
\]
given as the invariant part of $N_\delta$ with respect to $F$.

\subsection{The Conley index of $\isoS$.}
\begin{figure}[tb]
  \includegraphics[width=0.60\textwidth]{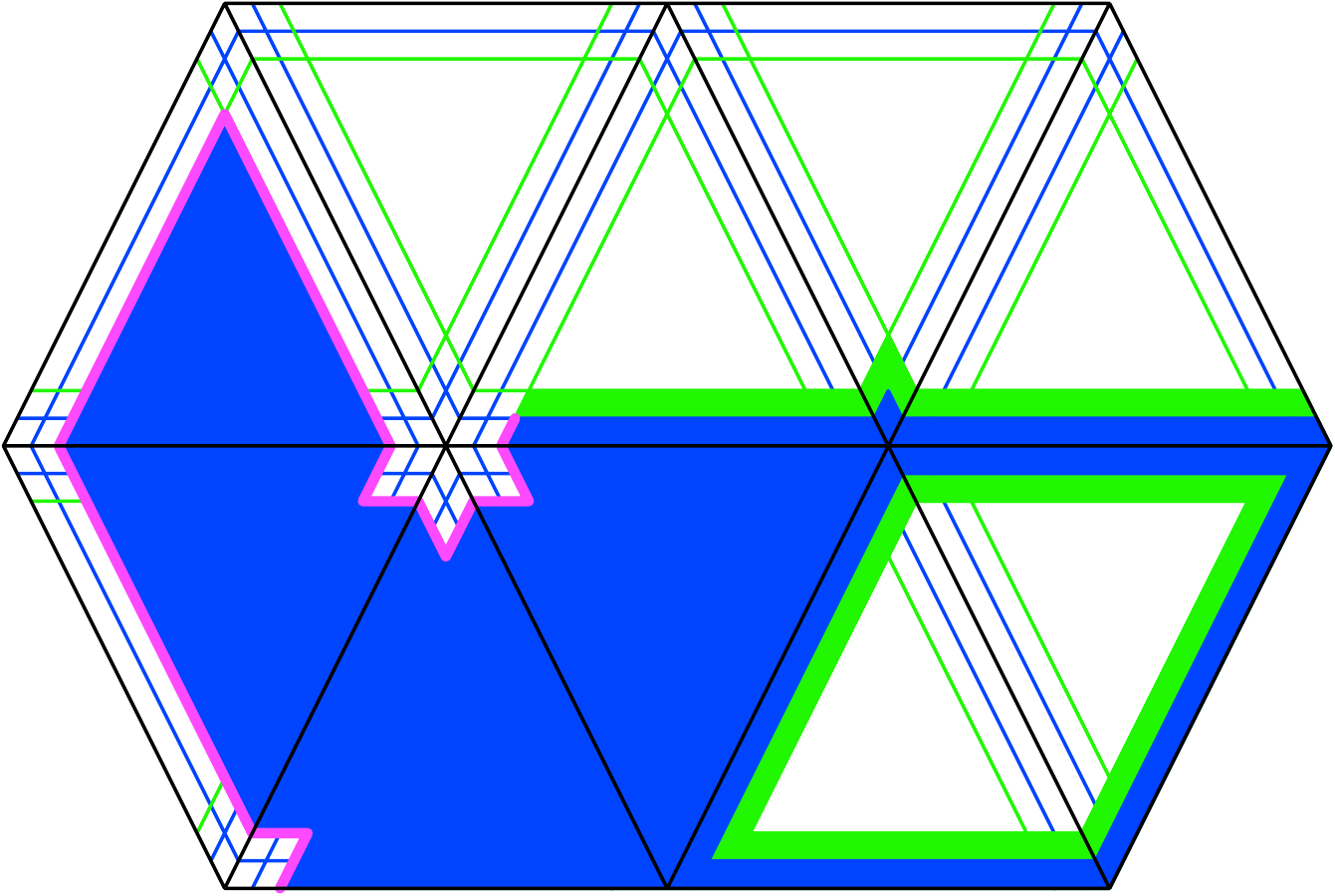}
  \caption{The weak index pair~$P = (P_1,P_2)$ associated with the
           isolating block~$N_\delta$ from Figure~\ref{fig:isoblock}.
           The set~$P_1$ is shown in dark blue, and the part
           of its boundary which comprises~$P_2$ is indicated
           in magenta. Notice that the parts of the $\delta$-cells
           shown in green are cut from the isolating block~$N_\delta$
           when passing to~$P_1$.}
  \label{fig:isopairP}
\end{figure}

In order to compare the Conley indices of $\cS$ and $\isoS$ we need to construct a weak index pair for $F$ in $N_\delta$.
To define such a weak index pair  we fix another constant $\delta'$ such that
\begin{equation} \label{def:deltagammaepsbound}
  0 < \delta' < \delta < \gamma < \epsilon < \frac{1}{d+1} \; ,
\end{equation}
and set
\begin{equation} \label{eq:P}
P_1:=N_\delta \cap N_{\delta'}
\qquad\mbox{ and }\qquad
P_2:= N_{\delta'} \cap \bd N_\delta\; .
\end{equation}
Clearly $P_2\subset P_1 \subset N:=N_\delta$ are compact sets.
We have the following theorem.
\begin{thm}\label{thm:P-is-weak-index-pair}
The pair $P=(P_1,P_2)$ defined by (\ref{eq:P}) is a weak index pair for $F$ and the isolating neighborhood $N=N_\delta$.
\end{thm}
The proof of Theorem~\ref{thm:P-is-weak-index-pair} will be presented in Section~\ref{P-is-weak-index-pair}.
A weak index pair for the  isolating block given in Figure~\ref{fig:isoblock}
is presented in Figure~\ref{fig:isopairP}.
As recalled in Section~\ref{sec:class-dyn}  the  Conley index of $\isoS$ with respect to  $F$ is
\[
\ClasCon(\isoS,F):=L(H^*(P),I_P),
\]
where $L$ is the Leray reduction of the relative cohomology graded module
$H^*(P)=H^*(P_1,P_2)$ of~$P$, and~$I_P$ is the index map on~$H^*(P)$.
In Section~\ref{sec:F-has-flowlike-Conley-index} we prove the following theorem.
\begin{thm}\label{thm:F-has-flowlike-Conley-index}
We have
\[
\ClasCon(\isoS,F)\cong (H^*(P),\id_{H^*(P)}),
\]
where $\id_{H^*(P)}$ denotes the identity map.
In other words, as in the case of flows, the Conley index of $\isoS$ with respect to $F$
can be simply defined as the relative cohomology $H^*(P)$.
\end{thm}

\subsection{Correspondence of Conley indices.}
As recalled in Section~\ref{sec:comb-dyn} the  Conley index of $\cS$ with respect to  $\Pi_\cV$ is
\[
\Con(\cS):=H^*(\ccl\cS,\Exit \cS)\; .
\]
In Section~\ref{sec:ConleyF=ConleyV} we prove the following theorem.
\begin{thm}\label{thm:ConleyF=ConleyV}
We have
\[
H^*(P_1,P_2) \cong H^*(\ccl\cS,\Exit \cS)\; .
\]
As a consequence,
\[
\ClasCon(\isoS) \cong \Con(\cS)\; .
\]
\end{thm}
Theorem~\ref{thm:ConleyF=ConleyV} extends the correspondence between the isolated invariant sets $\cS$ and $\isoS$
to the respective Conley indices.

\subsection{Correspondence of Morse decompositions.}

Given $\cM=\{\cM_r\,|\,r\in\PP\}$, a Morse decomposition of~$\cX$ with respect to the combinatorial flow~$\Pi_\cV$,
we define the sets
\[
M_r:=N^r _\epsilon\cap\ogr{\cM_r} \; ,
\]
where $N^r _\epsilon:=N_\epsilon(\cM_r)$ is given by (\ref{eq:N-beta}), that is, we have
\[
N^r _\epsilon=\bigcup_{\sigma\in\cM_r}\cl\cse{\sigma} \; .
\]
In Section~\ref{sec:Mcomb=MF} we prove the following theorem, which establishes the correspondence between
Morse decompositions of~$\cV$ and of~$F$.
\begin{thm}\label{thm:Mcomb=MF}
The collection $M:=\{M_r\,|\,r\in\PP\}$ is a Morse decomposition of $X$ with respect to $F$.
Moreover, for each $r\in\PP$ we have
$$
\Con (\cM_r)=C (M_r) \; ,
$$
and the Conley-Morse graphs for the Morse  decompositions $\cM$ and $M$ coincide.
\end{thm}

The reader can immediately see that Theorem~\ref{thm:main}, the main result of the
paper, is now an easy consequence of Theorems~\ref{thm:ConleyF=ConleyV}
and~\ref{thm:Mcomb=MF}.

\section{Proof of Theorem~\ref{thm:P-is-weak-index-pair}}
\label{P-is-weak-index-pair}

In this section we prove Theorem~\ref{thm:P-is-weak-index-pair}.
The proof is split into six auxiliary lemmas and the verification
that the pair~$P$ defined by (\ref{eq:P})
satisfies the conditions~(a) through~(d) of Definition~\ref{def:weak-index-pair}.

\subsection{Auxiliary lemmas.}
\begin{lem}\label{lem:A}
Consider a $\tilde{\delta}$ satisfying $0<\tilde{\delta}<\gamma$ and assume $A_\sigma \cap \cl\cstd{\tau}\neq \emptyset$ for all simplices $\tau,\sigma\in \cX$.
Then either $\tau$ is a face of $\sigma^-$ or $\tau=\sigma^+$.
\end{lem}
\proof Choose an $x\in A_\sigma \cap \cl\cstd{\tau}$. Accordingly to \eqref{def:ABCsigma} we have
\[
A_\sigma = \left\{ \tilde{x} \in \sigma^+\; \mid \;
    t_v(\tilde{x}) \geq \gamma \;\mbox{ for all }\; v \in \sigma^- \right\}
    \; \cup \; \sigma^-\; .
\]%
\begin{figure}[tb]%
  \includegraphics[width=0.40\textwidth]{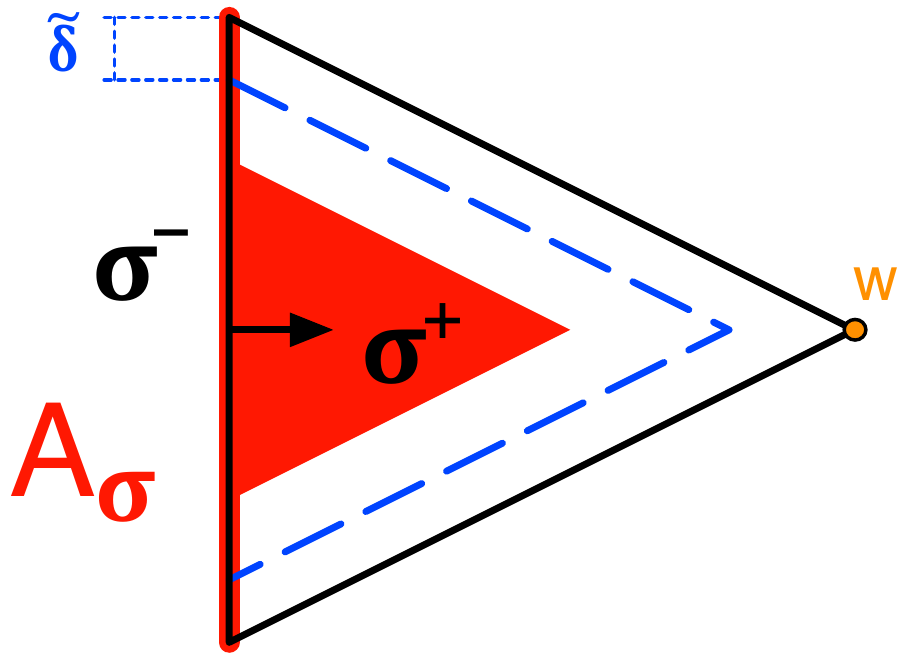}
  \caption{The sets $\sigma^+$, $A_\sigma$ and vertex $w$ in the proof of Lemma~\ref{lem:A}.}
  \label{fig:Asigmaw}
\end{figure}%
If $\tau$ is a face of $\sigma^-$, we are done. Suppose that this does not hold. Then $\tau$ has to contain the vertex $w$ of $\sigma^+$ complementing $\sigma^-$ as shown in Figure~\ref{fig:Asigmaw} and $x\notin \sigma^-$. This implies that $t_v(x)\geq \gamma > \tilde{\delta}$ for all $v\in\sigma^-$. Since
\[
\cl\cstd{\tau} = \left\{\tilde{x}\in X \; \mid \;
    t_v(\tilde{x}) \ge \tilde{\delta} \;\mbox{ for all }\; v \in \tau
    \;\mbox{ and }\;
    t_v(x) \le \tilde{\delta}
    \;\mbox{ for all }\; v \notin \tau \right\},
\]
this implies that all vertices of $\sigma^-$ have to be in $\tau$. Hence $\tau=\sigma^+$ and the claim is proved.
\qed
\begin{lem}\label{lem:B}
Suppose that $x\in \cl\csd{\tau}$ for some $\tau\in \cS$ and that $\sigma:=\sigma^\epsilon_{max}(x) \notin \cS$.
Then for every $\tilde{\delta}$ satisfying $0<\tilde{\delta} < \gamma$ we have
\[
F(x)\cap N_{\tilde{\delta}}=\emptyset \; .
\]
\end{lem}
\proof
Suppose that $F(x)\cap N_{\tilde{\delta}}\neq \emptyset$.
Hence, there exists a simplex $\hat{\tau}\in\cS$ and a point $y\in F(x)\cap \cl\cstd{\hat{\tau}}$. Then Lemma~\ref{lem:epsilon-delta-include} and Corollary~\ref{cor:equiv-cell-closure} imply that
\[
\sigma = \sigma^\epsilon_{max}(x) \subset  \sigma^\delta_{min}(x) \subset \tau.
\]
In other words, $\sigma$ is a face of $\tau$. Since $\cS$ is an isolated invariant set, $\tau\in\cS$ implies the inclusions $\tau^{\pm}\in \cS$.

Since $y\in F(x)$, we have $y\in F_\rho(x)$ for some simplex $\rho\in \cX^\epsilon(x)$. There are four possible cases to consider.

First, assume that $\rho\neq \sigma^+$ and $\rho\neq \sigma^-$. Then $F_\rho(x)=A_\rho$ and
$y\in A_\rho \cap \cl\cstd{\hat{\tau}}$. According to Lemma~\ref{lem:A}, $\hat{\tau}$ has to be a face of $\rho^-$ or $\hat{\tau}=\rho^+$. Recall that we assumed $\sigma\notin \cS$. Since $\rho\in\cX^\epsilon(x)$, we get
\[
\rho\subset \sigma \subset \tau.
\]
Since $\sigma\notin \cS$ and $\tau\in\cS$, Corollary~\ref{cor:convexity} implies that $\rho\notin \cS$
and Proposition~\ref{prop:exSclosed+-} implies that $\rho^{\pm}\notin\cS$.
Since we assumed $\hat{\tau}\in\cS$, Lemma~\ref{lem:A} shows that we cannot have $\hat{\tau}=\rho^+$. Thus, $\hat{\tau}$ has to be a face of $\rho^-$. The inclusions
\[
\hat{\tau} \subset \rho^- \subset \rho \subset \sigma \subset \tau,
\]
with
$\hat{\tau}, \tau \in \cS$ and $\sigma\notin \cS$ contradict the closedness of $\Exit\cS$ in Definition~\ref{def:discrete-isolated-inv}.

Now assume that $\rho = \sigma^+ \neq \sigma^-$. Since $\rho\in \cX^\epsilon(x)$, $\rho$ has to be a face of $\sigma$, so we get $\sigma^+=\sigma \notin \cS$ as well as $\sigma^-\notin \cS$. Moreover, in this case,
\[
y\in F_\rho(x)=B_\rho = B_{\sigma^+} \subset \sigma^+=\sigma.
\]
Hence, $\sigma^0_{min}(y)\subset \sigma$. Given that $y\in \cl\cstd{\hat{\tau}}$, we obtain
\[
\hat{\tau} \subset \sigma^{\tilde{\delta}}_{max}(y)\subset \sigma^0_{min}(y) \subset \sigma \subset \tau\; .
\]
This contradicts Corollary~\ref{cor:convexity}, because $\hat{\tau}, \tau \in \cS$ and $\sigma\notin \cS$.

Next assume that $\rho = \sigma^- \neq \sigma^+$. Then $y\in F_\rho(x)=C_\rho\subset A_\rho$. Hence, the inclusion $y\in A_\rho \cap \cl\cstd{\hat{\tau}}$ holds, and we get a contradiction as in the first case.

The last possible case is $\rho = \sigma^- =\sigma \in \Fix \cV$. Then $y \in \rho$. Since,
\[
\hat{\tau}\subset \sigma^0_{min}(y)\subset \sigma,
\]
we get the inclusions
\[
\hat{\tau} \subset \rho \subset \sigma \subset \tau,
\]
and a contradiction is reached as before.
\qed

\begin{lem} \label{lem:FdefTepsilon} (see \cite[Lemma 4.10]{KaMrWa16})
The image~$F(x)$ can be expressed alternatively as
\begin{equation} \label{eq:defF-2}
  F(x) \; = \;
  F_{\sigma_{max}^\epsilon(x)} \cup
    \bigcup_{\tau \in \cT^\epsilon(x)} F_\tau(x) \; ,
\end{equation}
where
\begin{equation}
\label{eq:cTeps}
  \cT^\epsilon(x) :=
  \left\{ \tau \in \cX^\epsilon(x) \setminus
    \left\{ \sigma_{max}^\epsilon(x) \right\} \; \mid \;
    \tau = \tau^-
    \;\mbox{ and }\;
    \tau^+ \notin \ccl \sigma_{max}^\epsilon(x)
  \right\} \; .
\end{equation}
Furthermore, every $\tau \in \cT^\epsilon(x)$ automatically
satisfies $\tau \in \dom\cV \setminus \Fix\cV$.
\end{lem}

\begin{lem}\label{lem:C}
We have
\[
F(N_\delta)\cap N_\delta \subset \ogr{\cS} .
\]
Consequently, $\Inv N_\delta \subset \ogr{\cS}$.
\end{lem}
\proof
We begin by showing that
\begin{equation}
\label{eq:C1}
   F(N_\delta)\cap N_\delta \subset |\cS| .
\end{equation}
Let $y\in F(N_\delta)\cap N_\delta$ with $y\in F(x)$ for some $x\in N_\delta$.
Then there exists a $\tau\in \cS$ and a $\hat{\tau}\in \cS$ such that $x\in \cl\csd{\tau}$ and $y \in \cl\csd{\hat{\tau}}$. Let $\sigma_0:=\sigma^0_{min}(x)$ and
$\sigma:=\sigma^\epsilon_{max}(x)$.
By Lemma~\ref{lem:epsilon-delta-include} and by Corollary~\ref{cor:equiv-cell-closure} we have
\[
\tau \subset \sigma^\delta_{max}(x) \subset \sigma^0_{min}(x) = \sigma_0,
\]
\[
\sigma = \sigma^\epsilon_{max}(x) \subset \sigma^0_{min}(x) = \sigma_0,
\]
and
\[
\sigma^\epsilon_{max}(x) \subset \sigma^\delta_{min}(x) \subset \tau\; .
\]
This implies that $\sigma\subset \tau \subset \sigma_0$.
We have $F(x)\cap N_\delta\neq \emptyset$, because $y\in F(x)\cap N_\delta$.
From  Lemma~\ref{lem:B} we get $\sigma\in \cS$.
Hence $\sigma^{\pm}\in\cS$.
Since $y\in F(x)$, we now consider the following cases resulting from Lemma~\ref{lem:FdefTepsilon}.

Assume first that  $y\in F_\rho(x)$ for a simplex $\rho \in \cT^\epsilon(x)$. Then $\rho\neq \sigma$, $\rho=\rho^-$, and $\rho^+\notin \ccl \sigma$. We will show that $\rho\in\cS$.
To see this assume the contrary.
Then one obtains $\rho \neq \sigma^-$, $\rho \neq \sigma^+$, so $F_\rho(x)=A_\rho$ and $y\in A_\rho\cap \cl\csd{\hat{\tau}}$ with $\hat{\tau}\in\cS$. Lemma~\ref{lem:A} implies that $\hat{\tau}=\rho^+$ or $\hat{\tau}$ is a face of $\rho$. If $\hat{\tau}=\rho^+$, then $\rho^+\in \cS$, which implies that  $\rho\in \cS$, a contradiction. If $\hat{\tau}$ is a face of $\rho$, we get the inclusions
\[
\hat{\tau} \subset \rho \subset \sigma^\epsilon_{max}(x)=\sigma\in\cS
\]
with $\hat{\tau},\sigma\in\cS$ and $\rho\notin\cS$. As in the proof of Lemma~\ref{lem:B}, this contradicts that $\cS$ is an isolated invariant set.
Thus,  $\rho\in\cS$. Then $\rho^+\in\cS$ and $y\in F_\rho(x) \subset \rho^+\subset|\cS|$.

Now assume that $y\in F_\sigma(x)=F_{\sigma^\epsilon_{max}(x)}(x)$. Then $F_\sigma(x)$ can either be~$B_\sigma$, or~$C_\sigma$, or~$\sigma$. All these sets are contained in $\sigma^+\in\cS$,
hence also in this case $y\in|\cS|$ and~\eqref{eq:C1} is proved.
By Proposition~\ref{prop:ogr}(ii), in order to conclude the proof it suffices to show that
\[
     F(N_\delta)\cap N_\delta\cap |\Exit\cS|=\emptyset.
\]
Assume the contrary. Then there exists a point $z\in F(N_\delta)\cap N_\delta$ with $z\in\cell{\tilde{\sigma}}$ for some $\tilde{\sigma}\in \Exit\cS$.
Since $z\in N_\delta$, there exists a $\tilde{\tau}\in \cS$ such that $z \in \cl \csd{\tilde{\tau}}$.
We get the inclusions
\[
\tilde{\tau} \subset \sigma^\delta_{max}(z) \subset \sigma^0_{min}(z)=\tilde{\sigma},
\]
with $\tilde{\tau}\in\cS$. Since $\tilde{\sigma}\in\Exit\cS$, this contradicts the closedness of $\Exit\cS$
and completes the proof.
\qed

\begin{lem} \label{lem:N-bd-cK}
Assume ~$\cS$ is an isolated invariant set for~$\Pi_\cV$ in the sense of
Definition~\ref{def:discrete-isolated-inv},
and consider the set $N=N_\delta\subset X = |\cX|$ given by \eqref{eq:N-delta}.
Then  $x \in \bd N_\delta$ if and only if
\begin{equation}
\label{eq:XdeltaS}
  \cX^\delta(x) \cap \cS \neq \emptyset
  \qquad\mbox{ and }\qquad
  \cX^\delta(x) \setminus \cS \neq \emptyset \; .
\end{equation}
\end{lem}
\proof
The fact that $x \in \bd N_\delta$ implies~\eqref{eq:XdeltaS} is shown
in~\cite[Lemma~5.5]{KaMrWa16}.
The reverse implication is an easy consequence of Proposition~\ref{prop:sigma-decomp}.
\qed

\begin{lem}\label{lem:D}
For any $x\in N_{\delta'} \cap \bd N_\delta$, we have $\sigma^\epsilon_{max}(x)\notin \cS$.
\end{lem}
\proof
Since $x \in \bd N_\delta$, Lemma~\ref{lem:N-bd-cK} implies that there exists a
$\sigma_1\in \cX^\delta(x)\setminus \cS$. Moreover, let $\sigma_0=\sigma^0_{min}(x)$ and $\sigma=\sigma^\epsilon_{max}(x)$.
By Corollary~\ref{cor:equiv-cell-closure} we then obtain the inclusions
$
\sigma=\sigma^\epsilon_{max}(x) \subset \sigma^\delta_{min}(x) \subset \sigma_1 \subset \sigma^\delta_{max}(x)\; .
$
In addition, we further have $x\in N_{\delta'}$.
Hence, there exists a simplex $\tau'\in\cS$ with $x\in\cl\csdp{\tau'}$. This implies that
$
\sigma^\delta_{max}(x) \subset \sigma^{\delta'}_{min}(x) \subset \tau' \subset \sigma^{\delta'}_{min}(x).
$
Together, this gives $ \sigma\subset \sigma_1 \subset \tau'$, where both $\sigma_1\notin\cS$ and $\tau'\in \cS$ hold.
Therefore, Corollary~\ref{cor:convexity}
 implies that $\sigma\notin \cS$, which completes the proof.
\qed

\subsection{Property (a).}
\begin{lem}\label{lem:E}
The sets $P_1,P_2$ given by \eqref{eq:P} and $N=N_\delta$ given by \eqref{eq:N-delta}
satisfy property (a) in Definition~\ref{def:weak-index-pair}.
\end{lem}
\proof
Let $i=1$. We know from Lemma~\ref{lem:C} that
\begin{equation}\label{eq:FP1capNsubS}
F(P_1)\cap N_\delta \subset F(N_\delta)\cap N_\delta \subset \ogr{\cS},
\end{equation}
and obviously $F(P_1)\cap N_\delta \subset N_\delta$. We argue by contradiction. Suppose that there exists an $x\in F(P_1)\cap N_\delta$ with $x\notin N_{\delta'}$.

Since $x \in N_\delta$, we have $x\in \cl\csd{\tau}$ for some $\tau \in \cS$.
Moreover, since $x \notin N_{\delta'}$, we have $x \in \cl\csdp{\tau'}$ for some $\tau' \notin \cS$.
Now let $\sigma_0=\sigma^0_{min}(x)$. Then one obtains the inclusions
$
\tau \subset \sigma^\delta_{max}(x) \subset \sigma^{\delta'}_{min}(x)\subset \tau' \subset \sigma^{\delta'}_{max}(x) \subset \sigma_0,
$
where $\tau \in \cS$ and $\tau' \notin \cS$.
Therefore, we get from Corollary~\ref{cor:convexity} that  $\sigma_0 \notin \cS$.
This shows that $x \notin \ogr{\cS}$, because by \eqref{eq:min-cell-sigma} we also have $x \in \cell\sigma_0$.
Since $x\in F(P_1)\cap N_\delta$, this contradicts (\ref{eq:FP1capNsubS}) and proves the claim for $i=1$.

Consider now the case $i=2$, and let $x\in P_2 = N_{\delta'}\cap \bd N_\delta$. Then Lemma~\ref{lem:D} implies $\sigma^\epsilon_{max}(x)\notin \cS$. Since $x\in N_\delta$, Lemma~\ref{lem:B} shows that
$
F(x)\cap N_\delta = \emptyset \subset P_2,
$
and the conclusion follows.
\qed
\subsection{Property (d).}
\begin{lem}\label{lem:F}
We have
\[
P_1 \setminus P_2 = N_{\delta'}\cap \Int N_\delta\; .
\]
As a consequence, property (d) in Definition~\ref{def:weak-index-pair} is satisfied.
\end{lem}
\proof
Let $x \in P_1 \setminus P_2$ be arbitrary. Since $x\in P_1$, we have $x \in N_\delta$ and $x \in N_{\delta'}$. Since $x \notin P_2$, either $x \notin \bd N_\delta$ or $x \notin N_{\delta'}$. The second case is excluded, hence we have $x \notin \bd N_\delta$. It follows that $x\in \Int N_\delta$, and therefore also that $x \in N_{\delta'}\cap \Int N_\delta$.

Conversely, let $x \in N_{\delta'}\cap \Int N_\delta \subset P_1$. Then both $x \notin \bd N_\delta$
and $x\notin P_2$ are satisfied. It follows that $x\in  P_1 \setminus P_2$.

Now, property (d) trivially follows from the inclusion $N_{\delta'}\cap \Int N_\delta \subset \Int N_\delta$.
\qed

\subsection{Property (c).}
\begin{lem}\label{lem:G}
We have
\[
\Inv N_\delta \subset \Int(P_1\setminus P_2).
\]
In other words, property (c) in Definition~\ref{def:weak-index-pair} is satisfied.
\end{lem}
\proof
According to Theorem~\ref{th:N-iso-block} we have
\begin{equation}\label{eq:Ndelta-iso-nb}
\Inv N_\delta \subset \Int N_\delta \; .
\end{equation}
We will show in the following that also
\begin{equation}\label{eq:Ndelta-Nprime}
\Inv N_\delta \subset \Int N_{\delta'} \; .
\end{equation}
We argue by contradiction. Suppose that $\Inv N_\delta \setminus \Int N_{\delta'}\neq \emptyset$ and let $x\in \Inv N_\delta$ be such that $x\notin \Int N_{\delta'}$.
If $x \in \bd N_{\delta'}$, then Lemma~\ref{lem:N-bd-cK} shows that there exists a simplex $\tau'\notin \cS$ with $x \in \cl\csdp{\tau'}$. It is clear that if $x \notin N_{\delta'}$ such a $\tau'$ also exists.
According to Lemma~\ref{lem:C}, we have $x \in \ogr{\cS}$, that is, there exists a simplex $\sigma\in \cS$ with $x \in \cell\sigma$.
By \eqref{eq:min-cell-sigma} we have $\sigma^0_{min}(x)=\sigma$. As before, we obtain inclusions
\[
\sigma^\epsilon_{max}(x) \subset \sigma^{\delta'}_{min}(x) \subset \tau' \subset \sigma^{\delta'}_{max}(x) \subset \sigma^0_{min}(x) = \sigma,
\]
with $\tau' \notin \cS$ and $\sigma \in \cS$, and Corollary~\ref{cor:convexity} implies
that $\sigma^\epsilon_{max}(x) \notin \cS$. Due to our assumption, we have $x \in \Inv N_\delta \subset N_\delta$, and Lemma~\ref{lem:B} gives $F(x)\cap N_\delta=\emptyset$. Hence $x \notin \Inv N_\delta$, which is a contradiction and thus proves \eqref{eq:Ndelta-Nprime}.

The inclusions (\ref{eq:Ndelta-iso-nb}) and (\ref{eq:Ndelta-Nprime}) give
\[
\Inv N_\delta \subset  \Int N_{\delta'} \cap \Int N_\delta \subset \Int (N_{\delta'}\cap \Int N_\delta) = \Int(P_1\setminus P_2),
\]
which completes the proof.
\qed

For the next result we need the following two simple observations.
If~$A$ and~$B$ are closed subsets of~$X$, then
\begin{equation}\label{eq:bdAandB}
\bd(A\cap B) \subset (\bd A \cap B) \cup (A \cap \bd B) \;.
\end{equation}
If $A$ is a closed subset of $X$, then
\begin{equation}\label{eq:bfFAsubAcapFA}
\bd_F(A) \subset A \cap F(A)\; .
\end{equation}
The first observation is straightforward. In order to verify the second one, it is clear that $\bd_F(A) \subset \cl A = A$. Since $F(A) \setminus A \subset F(A)$, we get
\[
\cl(F(A)\setminus A)\subset \cl F(A) = F(A),
\]
because the map~$F$ is upper semi-continuous and the set~$X$ is compact. These two inclusions immediately give~(\ref{eq:bfFAsubAcapFA}).

\subsection{Property (b).}
\begin{lem}\label{lem:H}
For $P_1$ and $P_2$ defined by (\ref{eq:P}), we have
\[
\bd_F(P_1)\subset P_2\; .
\]
In other words, property (b) in Definition~\ref{def:weak-index-pair} is satisfied.
\end{lem}
\proof
One can easily see that $\bd_F(P_1) \subset \bd(P_1)$. Together with (\ref{eq:bdAandB})
this further implies
\begin{eqnarray*}
\bd_F(P_1) & \subset & \bd(P_1) = \bd(N_\delta \cap N_{\delta'}) \\
           & \subset & (N_{\delta'} \cap \bd N_\delta)\cup (N_\delta \cap \bd N_{\delta'})\\
           & = & P_2 \cup (N_\delta \cap \bd N_{\delta'})\; .
\end{eqnarray*}
Thus, if we can show that
\begin{equation}\label{eq:bdFP1capempty}
\bd_F(P_1) \cap (N_\delta \cap \bd N_{\delta'}) = \emptyset,
\end{equation}
then the proof is complete. We prove this by contradiction. Assume that there exists an $x \in \bd_F(P_1) \cap (N_\delta \cap \bd N_{\delta'})$.
Since $x \in \bd_F(P_1)$, by (\ref{eq:bfFAsubAcapFA}) we get
\[
x \in P_1 \cap F(P_1) \subset N_\delta \cap F(N_\delta)\; .
\]
Thus, due to Lemma~\ref{lem:C}, we have the inclusion $x \in \ogr{\cS}$.
It follows that  there exists a simplex $\sigma \in \cS$ with $x \in \cell\sigma$ and by \eqref{eq:min-cell-sigma}
$\sigma=\sigma^0_{min}(x)$.

Since $x\in N_\delta$, we  get a simplex $\tau \in \cS$ such that $x \in \cl\csd{\tau}$
and since $x \in \bd N_{\delta'}$, by Lemma~\ref{lem:N-bd-cK} we also get a simplex $\tau' \notin \cS$ such that $x \in \cl\csdp{\tau'}$. Now,  Lemma~\ref{lem:epsilon-delta-include} and Corollary~\ref{cor:equiv-cell-closure}
imply
\[
\tau \subset \sigma^\delta_{max}(x) \subset \sigma^{\delta'}_{min}(x) \subset \tau' \subset  \sigma^{\delta'}_{max}(x) \subset \sigma^0_{min}(x) = \sigma \; .
\]
Since~$\tau, \sigma \in \cS$~ and $\tau' \notin \cS$,
this contradicts, in combination with Corollary~\ref{cor:convexity},
the fact that~$\cS$ is an isolated invariant set --- and the proof is complete.
\qed

\medskip

Theorem~\ref{thm:P-is-weak-index-pair} is now an immediate consequence of
Lemmas \ref{lem:E}, \ref{lem:F}, \ref{lem:G}, \ref{lem:H}.

\section{Proof of Theorem~\ref{thm:F-has-flowlike-Conley-index}}
\label{sec:F-has-flowlike-Conley-index}

In this section we prove Theorem~\ref{thm:F-has-flowlike-Conley-index}.
Since the Leray reduction of an identity is clearly the same identity, it suffices
to prove that the index map~$I_P$ is the identity map. We achieve this by constructing
an acyclic-valued and upper semicontinuous map~$G$ whose graph contains both the graph
of~$F$ and the graph of the identity. The map~$G$ is constructed by gluing
two multivalued and acyclic maps. One of these, the map~$\tilde{F}$ defined below, is a modification of our map~$F$, while the second map~$D$ contains
the identity.

\subsection{The map $\tilde{F}$.}
For $x\in X$ and  $\sigma\in\cX^\epsilon(x)$ we define
\[
  \tilde{F}_\sigma(x):=F_\sigma(x)\cup A_\sigma.
\]
One can then easily verify that
\begin{equation}
\label{eq:Ftilde-formula-zero}
\tilde{F}_\sigma(x)=\left\{\begin{array}{ll}
A_\sigma&\mbox{ if }F_\sigma(x)=C_\sigma\\
\sigma^+&\mbox{ if }F_\sigma(x)=B_\sigma\\
F_\sigma(x)&\mbox{ otherwise},
\end{array}
\right.
\end{equation}
and that the inclusion
\begin{equation}\label{eq:sigma_subset_Ftilde_sigma}
\sigma\subset\tilde{F}_\sigma(x)
\end{equation}
holds.
We will show that the auxiliary map $\tilde{F}:X\mto X$ given by
\begin{equation}
\label{eq:Ftilde}
\tilde{F}(x):=\bigcup_{\sigma\in\cX^\epsilon(x)}\tilde{F}_\sigma(x).
\end{equation}
is acyclic-valued. For this we need a few auxiliary results.
\begin{lem}\label{lem:Ftilde-formula-one}
  For any $x\in X$ and $\sigma:=\sigma^{\epsilon}_{max}(x)$ we have
\begin{equation}
\label{eq:Ftilde-formula-one}
\tilde{F}(x)=A_\sigma\cup F(x).
\end{equation}
\end{lem}
\proof
  It is straightforward to observe that the right-hand side of \eqref{eq:Ftilde-formula-one}
is contained in the left-hand side. To prove the opposite inclusion, take a $y\in \tilde{F}(x)$
and select a simplex $\tau\in\cX^\epsilon(x)$ such that $y\in \tilde{F}_\tau(x)$.
In particular, $\tau\subset\sigma$.
Note that if $\tau=\sigma$, then
\begin{equation}
\label{eq:Ftilde-formula-one-1}
\tilde{F}_\tau(x)=\tilde{F}_\sigma(x)=F_\sigma(x)\cup A_\sigma\subset F(x)\cup A_\sigma.
\end{equation}
According to \eqref{def:Fsigma} the value $F_\tau(x)$ may be $A_\tau$, $B_\tau$, $C_\tau$ or $\tau$.
Assume first that we have $F_\tau(x)=A_\tau$. Then,
\[
  \tilde{F}_\tau(x)=F_\tau(x)\cup A_\tau=F_\tau(x)\subset F(x)\subset F(x)\cup A_\sigma.
\]
Next, consider the case $F_\tau(x)=B_\tau$. Then, $\tau=\sigma^+\neq\sigma^-$.
Since $\tau\subset\sigma$, we cannot have $\sigma=\sigma^-$. Hence, $\sigma=\sigma^+=\tau$
and estimation \eqref{eq:Ftilde-formula-one-1} applies.
Assume in turn that $F_\tau(x)=C_\tau$. Then, $\tau=\sigma^-\neq\sigma^+$,  $A_\tau=A_\sigma$
and
\[
  \tilde{F}_\tau(x)=C_\tau\cup A_\tau=A_\tau=A_\sigma\subset F(x)\cup A_\sigma.
\]
Finally, if $F_\tau(x)=\tau$, then $\tau=\sigma^-=\sigma^+=\sigma$
and again \eqref{eq:Ftilde-formula-one-1} applies.
\qed

\begin{prop}
\label{prop:tau-sigma}
  Let $x\in X$ and let $\sigma:=\sigma^{\epsilon}_{max}(x)$.
 For any simplex $\tau\in\cT^\epsilon(x)$, where~$\cT^\epsilon(x)$ is given
 as in~\eqref{eq:cTeps}, we have
\begin{equation}
\label{eq:Ftilde-formula-two-2}
 A_\tau\cap A_\sigma=\tau\cap \sigma^-.
\end{equation}
\end{prop}
\proof
We have $\tau\cap\sigma^-\subset\sigma^-\subset A_\sigma$ and $\tau\cap\sigma^-\subset\tau=\tau^-\subset A_\tau$,
which shows that the right-hand side of \eqref{eq:Ftilde-formula-two-2} is contained in the left-hand side.
Observe that for any simplex $\rho\neq\sigma^+$ we have $\rho\cap A_\sigma\subset \sigma^-$.
We cannot have $\tau^+=\sigma^+$, because then either $\sigma=\sigma^-=\tau^-=\tau$ or $\tau^+=\sigma^+=\sigma\in\ccl\sigma$, in both cases contradicting the inclusion $\tau\in\cT^\epsilon(x)$.
Thus, $A_\tau\cap A_\sigma\subset \tau^+\cap A_\sigma\subset\sigma^-$.
It follows that we must have $A_\tau\cap A_\sigma\subset\sigma^-\cap\tau^+\subset\sigma\cap\tau^+$.
The simplex $\sigma\cap\tau^+$ must be a proper face of $\tau^+$, because otherwise $\tau^+$ is a face of
$\sigma$ which contradicts $\tau^+\not\in\ccl\sigma.$ But, $\tau\subset\sigma\cap\tau^+$ and $\tau=\tau^-$
is a face of $\tau^+$ of codimension one.
It follows that $\sigma\cap\tau^+=\tau$ which proves \eqref{eq:Ftilde-formula-two-2}.
\qed

\medskip

The following proposition is implicitly proved
in the second to last paragraph of the proof of \cite[Theorem 4.12]{KaMrWa16}.
\begin{prop}
\label{prop:Ftau-Atau}
   For any $\tau\in\cT^\epsilon(x)$ we have $F_\tau(x)=A_\tau$.
\qed
\end{prop}

\begin{lem}\label{lem:Ftilde-formula-two}
  For an  $x\in X$ and $\sigma:=\sigma^{\epsilon}_{max}(x)$ we have
\begin{equation}
\label{eq:Ftilde-formula-two}
A_\sigma\cap F(x)=A_\sigma\cap F_\sigma(x).
\end{equation}
\end{lem}
\proof
   Obviously, the right-hand side is contained in the left-hand side. To prove the opposite inclusion,
 choose a $y\in F(x)\cap A_\sigma$ and select a simplex~$\tau$ such that $y\in F_\tau(x)$.
 It suffices  to show that $y\in F_\sigma(x)$.
 By Lemma~\ref{lem:FdefTepsilon} we may assume that either $\tau=\sigma$ or $\tau\in\cT^\epsilon(x)$.
 If $\tau=\sigma$, the inclusion is obvious. Hence, assume that $\tau\in\cT^\epsilon(x)$.
 By Corollary~\ref{cor:equiv-cell-closure} we have $\tau\subset\sigma$.
 This means that $\tau\subsetneq\sigma$, $\tau=\tau^-$ and $\tau^+\not\in\ccl\sigma$.
 From Proposition~\ref{prop:Ftau-Atau} we get $F_\tau(x)=A_\tau$.
 Therefore,
\begin{equation}
\label{eq:Ftilde-formula-two-1}
 y\in A_\tau\cap A_\sigma\subset \tau^+\cap A_\sigma\subset\tau^+\cap \sigma^+.
\end{equation}
If $\sigma=\sigma^+=\sigma^-$, then $y\in A_\sigma=\sigma=F_\sigma(x)$, hence the inclusion holds.
Thus, consider the case $\sigma^+\neq \sigma^-$.
By Proposition~\ref{prop:tau-sigma} we get $y\in\tau\cap\sigma^-$.
We cannot have $\tau\cap\sigma^-=\sigma^-$, because then $\sigma^-\subset\tau\subsetneq\sigma$,
$\tau^-=\tau=\sigma^-$, $\tau^+=\sigma^+$ and $\tau\neq\sigma$ implies $\tau^+=\sigma^+=\sigma$,
$\tau^+\in\ccl\sigma$, a contradiction. Hence, $\tau\cap\sigma^-$ is a proper face of the simplex $\sigma^-$.
Therefore, $\tau\cap\sigma^-\subset C_\sigma\subset F_\sigma(x)$.
\qed

\begin{thm}
\label{thm:Ftilde}
  The map $\tilde{F}$ is upper semicontinuous and acyclic-valued.
\end{thm}
\proof
   The upper semicontinuity of the map~$\tilde{F}$ is an immediate consequence of
   formula~\eqref{eq:Ftilde} and Lemma~\ref{lem:cs-semicontiniuty}. To show that $\tilde{F}$
is acyclic-valued fix an $x\in X$. By \eqref{eq:Ftilde-formula-one},
$\tilde{F}(x)=A_\sigma\cup F(x)$. The set $A_\sigma$ is acyclic by Lemma~\ref{lem:ABC-contractible}
and the set $F(x)$ is acyclic by Theorem~\ref{th:F-usc-acyclic}. Moreover,
$A_\sigma\cap F(x)=A_\sigma\cap F_\sigma(x)$ by \eqref{eq:Ftilde-formula-two}.
Hence, due to~\eqref{def:Fsigma} the intersection
$A_\sigma\cap F(x)$ is either $A_\sigma$ or $C_\sigma$, hence also acyclic. Therefore, it follows
from the Mayer-Vietoris theorem that $\tilde{F}(x)$ is acyclic.
\qed

\begin{prop}\label{prop:tildeF_pos_inv}
The weak index pair~$P$ is positively invariant with respect to~$\tilde{F}$ and~$N_\delta$, that is, we have
$$
\tilde{F}(P_i)\cap N_\delta\subset P_i
\quad\mbox{ for }\quad
i=1,2.
$$
\end{prop}
\proof The proof is analogous to the proof of Lemma \ref{lem:E}.
\qed\\


\subsection{The map $\tilde{D}$.}
We define a multivalued map $D:X\mto X$, by letting
$$
D(x):=\conv(\{x\}\cup\sigma ^\varepsilon _{max}(x)),
$$
where $\conv A$ denotes the convex hull of $A$.
Note that the above definition is well-posed, because both $\{x\}$ and $\sigma$ are subsets of the same simplex $\sigma ^0 _{min}(x)$.

In order to show that $D$ is upper semicontinuous we need the following lemma.
\begin{lem}
\label{lem:sigma-max-usc}
  The mapping
\[
  X\ni x\mapsto\sigma^\epsilon_{max}(x)\subset X
\]
is strongly upper semicontinuous, that is, for every $x\in X$ there exists a neighborhood~$V$ of~$x$
such that for each $y\in V$ we have $\sigma^\epsilon_{max}(y)\subset \sigma^\epsilon_{max}(x)$.
\end{lem}
\proof
  By Lemma~\ref{lem:cs-semicontiniuty}, we can choose a neighborhood~$V$ of~$x$ in such a way that
  $\cX^\epsilon(y)\subset\cX^\epsilon(x)$ for $y\in V$. In particular,
  $\sigma^\epsilon_{max}(y)\in \cX^\epsilon(y)\subset\cX^\epsilon(x)$.
  By Corollary~\ref{cor:equiv-cell-closure}, we obtain $\sigma^\epsilon_{max}(y)\subset \sigma^\epsilon_{max}(x)$.
\qed
\begin{prop}\label{prop:D_usc_contractible}
The mapping $D$ is upper semicontinuous and has non-empty and contractible values.
\end{prop}
\proof
  Since the values of $D$ are convex, they are obviously contractible.
  To see that $D$ is upper semicontinuous, fix $\epsilon>0$.
  By Lemma~\ref{lem:sigma-max-usc} we can find a neighborhood $V$ of $x$ such that
  $\sigma^\epsilon_{max}(y)\subset \sigma^\epsilon_{max}(x)$.
  Let $B(x,\epsilon)$ denote the $\epsilon$-ball around $x$, let $y\in B(x,\epsilon)\cap V$,
  and fix a point $z\in D(y)$. Then, $z=ty+(1-t)\bar{y}$ for a $t\in[0,1]$ and $\bar{y}\in \sigma^\epsilon_{max}(y)$.
  Let $z':=tx+(1-t)\bar{y}$. Since $\sigma^\epsilon_{max}(y)\subset \sigma^\epsilon_{max}(x)$, we have $z'\in D(x)$.
  Moreover, the estimate $||z-z'||=t||y-x||\leq ||y-x||<\epsilon$ holds. It follows that $z\in B(D(x),\epsilon)$.
  Hence, $D(y)\subset B(D(x),\epsilon)$, which proves the upper semicontinuity of $D$.
\qed

\begin{lem}\label{lem:aux_D_inv}
Let $0<\tilde{\delta}<\varepsilon$. For any $x\in X$ and any $y\in D(x)$, $y\neq x$, we have
$$
\sigma ^{\tilde{\delta}} _{max}(y)\subset\sigma ^{\tilde{\delta}} _{min}(x).
$$
\end{lem}
\proof
Let $x\in X$ and $y\in D(x)$, with $y\neq x$, be fixed. Then $y=\alpha x+(1-\alpha)x_\sigma$ for some $x_\sigma\in\sigma^\epsilon_{max}(x)$ and $\alpha\in[0,1)$. Consider a vertex $v\notin\sigma ^{\tilde{\delta}} _{min}(x)$. By Lemma~\ref{lem:epsilon-delta-include} we have $\sigma ^\varepsilon _{max}(x)\subset\sigma ^{\tilde{\delta}} _{min}(x)$, which shows that $v\notin \sigma ^\varepsilon _{max}(x)$. Therefore,
$$
\begin{array}{rcl}
t_v(y)&=&\alpha t_v(x)+(1-\alpha)t_v(x_\sigma)\\[1ex]
&=&\alpha t_v(x)<t_v(x)\leq\tilde{\delta},
\end{array}
$$
which implies $v\notin\sigma ^{\tilde{\delta}} _{max}(y)$, and the inclusion $
\sigma ^{\tilde{\delta}} _{max}(y)\subset\sigma ^{\tilde{\delta}} _{min}(x)
$
follows.
\qed
\begin{prop}\label{prop:D_pos_inv}
The weak index pair~$P$ is positively invariant with respect to~$D$ and~$N_\delta$, that is,
we have
\begin{equation}
\label{eq:D_pos_inv}
D(P_i)\cap N_\delta\subset P_i
\quad\mbox{ for }\quad
i=1,2.
\end{equation}

\end{prop}
\proof
We begin with the proof for the case $i=1$. Fix an  $x\in P_1=N_\delta\cap N_{\delta'}$.
Then $x\in\cl\csd{\tau}$ and $x\in\cl\csdp{\tau'}$ for some $\tau, \tau'\in\cS$.
Consider a $y\in D(x)\cap N_\delta$. For $y=x$ inclusion \eqref{eq:D_pos_inv} trivially holds.
Therefore, assume $y\neq x$. Since $y\in N_\delta$, there exists a simplex $\eta\in\cS$ with $y\in\cl\csd{\eta}$. Consider a simplex $\eta'\in\cX$ such that $y\in\cl\csdp{\eta'}$. By Lemma \ref{lem:aux_D_inv} we obtain
$
\eta'\subset\sigma^{\delta'}_{max}(y)\subset\sigma^{\delta'}_{min}(x)\subset\tau'
$
and Lemma~\ref{lem:epsilon-delta-include} together with Corollary~\ref{cor:equiv-cell-closure} implies
$
\eta\subset\sigma^{\delta}_{max}(y)\subset\sigma^{\delta'}_{min}(y)\subset\eta'.
$
Thus,
$
\eta\subset\eta'\subset\tau'.
$
Since $\eta\in\cS$ and $\tau'\in\cS$, we get from Corollary~\ref{cor:convexity} that $\eta'\in \cS$.
Thus, $y\in N_{\delta'}$, which completes the proof for $i=1$.

To proceed with the proof for $i=2$, fix an $x\in P_2$ and consider a $y\in D(x)\cap N_\delta$.
Inclusion \eqref{eq:D_pos_inv} is trivial when $y=x$.
Hence, assume $y\neq x$. Note that $P_2\subset N_\delta'$, therefore $x\in\cl\csdp{\tau'}$ for some $\tau'\in\cS$. Since $x\in\bd N_\delta$, by Lemma~\ref{lem:N-bd-cK} there exists a simplex $\tau\notin \cS$ such that $x\in\cl\csd{\tau}$. According to Lemma~\ref{lem:epsilon-delta-include} we then have the inclusion
$
\tau\subset\sigma^{\delta}_{max}(x)\subset\sigma^{\delta'}_{min}(x)\subset\tau',
$
and this yields $\tau\in\Exit\cS$. Consider a simplex $\eta$ such that $y\in\cl\csd{\eta}$. By Lemma \ref{lem:aux_D_inv} we have
$
\eta\subset\sigma^{\delta}_{max}(y)\subset\sigma^{\delta}_{min}(x) \subset\tau.
$
Now, the closedness of $\Exit\cS$ implies $\eta\in\Exit\cS$.
Hence, $y\in\bd N_\delta$ by Lemma~\ref{lem:N-bd-cK}. Observe that by case $i=1$ we also have $y\in P_1\subset N_\delta'$. Therefore, $y\in P_2$.
\qed\\


\subsection{The map $G$.}
Define the multivalued  map $G:X\mto X$ by
\begin{equation}\label{eq:defG}
G(x):=D(x)\cup\tilde{F}(x).
\end{equation}
\begin{prop}\label{lem:D_properties}
The following conditions hold:
\begin{itemize}
\item[(i)] $G$ is upper semicontinuous,
\item[(ii)] $P$ is positively invariant with respect to $G$ and $N_\delta$,
\item[(iii)] $G$ is acyclic-valued.
\end{itemize}
\end{prop}
\proof
The map $G$ inherits properties (i) and (ii) directly from its summands $\tilde{F}$ and $D$ (see Theorem \ref{thm:Ftilde}, Proposition~\ref{prop:tildeF_pos_inv}, Proposition \ref{prop:D_usc_contractible} and Proposition \ref{prop:D_pos_inv}).

To prove (iii), fix an $x\in X$ and note that $\tilde{F}(x)$ is acyclic by Theorem~\ref{thm:Ftilde},
and that~$D(x)$ is acyclic by Proposition~\ref{lem:D_properties}(iii).
Let $\sigma:=\sigma^\epsilon_{max}(x)$ and $\sigma^0:=\sigma^0_{min}(x)$.
Obviously either $x\in\sigma$ or $x\not\in\sigma$.
To begin with, we consider the case $x\in\sigma$. Then one has $D(x)=\sigma$.
We will show that $D(x)\subset \tilde{F}(x)$. Indeed, if $\sigma=\sigma^-$ then we have $D(x)=\sigma=\sigma^-\subset A_\sigma\subset \tilde{F}(x)$ by Lemma \ref{lem:Ftilde-formula-one}.
If $\sigma\neq\sigma^-$, then one has the equality $\sigma=\sigma^+$.
In that case $F_\sigma(x)=B_\sigma$ and by \eqref{eq:Ftilde-formula-zero}
we get $\tilde{F}_\sigma(x)=\sigma^+$, which shows that $D(x)=\sigma=\sigma^+=\tilde{F}_\sigma(x)\subset
\tilde{F}(x)$. Consequently, if $x\in\sigma$, then we have the equality $G(x)=\tilde{F}(x)$,
and this set is acyclic by Theorem \ref{thm:Ftilde}.

Thus, consider the case $x\notin\sigma$.
By the Mayer-Vietoris theorem it suffices to show that $D(x)\cap \tilde{F}(x)$ is acyclic, because
both sets $D(x)$ and $\tilde{F}(x)$ are acyclic.
To this end we use the following representation
$$
\tilde{F}(x)=F_\sigma(x)\cup A_\sigma\cup\bigcup_{\tau\in\cT^\varepsilon(x)} F_\tau(x),
$$
which follows immediately from Lemma \ref{lem:Ftilde-formula-one} and \eqref{eq:defF-2}.
First, we will show that for every simplex $\tau\in\cT^\epsilon(x)$ we have
\begin{equation}
\label{eq:D_properties-1}
    D(x)\cap F_\tau(x) \subset D(x)\cap \left(F_\sigma(x)\cup A_\sigma\right).
\end{equation}
To see this observe that
since $x\in\sigma^0\setminus\sigma$, it is evident that $\sigma$ is a proper face of
the simplex~$\sigma^0$. Moreover, one can easily observe that
\begin{equation}\label{eq:D_int_bd_eta=sigma}
D(x)\subset\sigma^0\mbox{ and }D(x)\cap |\Bd\sigma^0|=\sigma.
\end{equation}
Note that from the definition of the collection $\cT^\varepsilon(x)$ (cf. Lemma~\ref{lem:FdefTepsilon}) it follows that any simplex $\tau\in\cT^\varepsilon(x)$  is a proper face of $\sigma$. Therefore, $\sigma^0$ is a coface of $\tau$ of codimension greater than one. Hence, we cannot have $\tau^+=\sigma^0$.
By Proposition~\ref{prop:Ftau-Atau} we have  $F_\tau(x)=A_\tau\subset\tau^+$.
Thus, from $\tau=\tau^-\subset A_\tau$ and (\ref{eq:D_int_bd_eta=sigma}), we obtain
\begin{equation}\label{eq:D_intersect_Ftau}
D(x)\cap F _\tau (x)\subset\sigma\mbox{ for any }\tau\in\cT ^\varepsilon(x).
\end{equation}
We will now show that
\begin{equation}
\label{eq:D_properties-2}
\sigma\subset D(x)\cap \left(F_\sigma(x)\cup A_\sigma\right).
\end{equation}
Obviously, $\sigma\subset D(x)$. If $\sigma=\sigma^-$, then
$\sigma\subset A_\sigma\subset F_\sigma(x)\cup A_\sigma$.
If $\sigma=\sigma^+$, then by~\eqref{def:Fsigma} we have $\sigma\subset\sigma^+\subset F_\sigma(x)\cup A_\sigma$.
Hence, \eqref{eq:D_properties-2} is proved.
Formula \eqref{eq:D_properties-1} follows now from \eqref{eq:D_intersect_Ftau} and \eqref{eq:D_properties-2}.
From \eqref{eq:D_properties-1} we immediately obtain that
\begin{equation}
\label{eq:D_properties-3}
D(x)\cap \tilde{F}_\tau (x)= D(x)\cap \left(F_\sigma(x)\cup A_\sigma\right).
\end{equation}
Now we distinguish the two complementary cases: $\sigma=\sigma^+$ and $\sigma=\sigma^-\neq\sigma^+$. First of all, if $\sigma=\sigma^+$, then
\[
   D(x)\cap \left(F_\sigma(x)\cup A_\sigma\right)\subset\sigma^0\cap\sigma^+=\sigma^0\cap\sigma=\sigma.
\]
It follows from \eqref{eq:D_properties-2} and \eqref{eq:D_properties-1} that
in this case
\[
    D(x)\cap\tilde{F}(x)=D(x)\cap \left(F_\sigma(x)\cup A_\sigma\right)=\sigma
\]
is an acyclic set. We show that the same is true in the second case $\sigma=\sigma^-\neq\sigma^+$.
Now one has $F_\sigma(x)\cup A_\sigma=C_\sigma\cup A_\sigma=A_\sigma$.
Observe that $A_\sigma=A^+\cup \sigma^-$, where
$$
A^+:=\{y\in\sigma ^+\;|\;t_v(y)\geq\gamma\mbox{ for }v\in\sigma^-\}
$$
is a convex set.
This, together with \eqref{eq:D_properties-3},
shows that
\begin{equation}
\label{eq:D_properties-4}
D(x)\cap \tilde{F}(x)=D(x)\cap A_\sigma=D(x)\cap (A^+\cup \sigma^-)=D(x)\cap A^+\cup \sigma^-.
\end{equation}
The acyclicity of the right-hand-side of \eqref{eq:D_properties-4} follows from the Mayer-Vietoris theorem,
because $D(x)\cap A^+$, $\sigma^-$, and $D(x)\cap A^+\cap\sigma^-=A^+\cap\sigma^-$ are all convex.
Therefore, by \eqref{eq:D_properties-4} also in this case the set $D(x)\cap \tilde{F}(x)$ is acyclic.
This completes the proof.
\qed\\
We are now able to prove Theorem~\ref{thm:F-has-flowlike-Conley-index}.

\medskip

{\bf Proof of Theorem~\ref{thm:F-has-flowlike-Conley-index}:}
By Lemma \ref{prop:tildeF_pos_inv} and Proposition \ref{prop:D_pos_inv} we can consider
the map~$G$, given by~(\ref{eq:defG}), as a map of pairs
$$
G:(P_1,P_2)\mto(T_1(P),T_2(P)).
$$
 Directly from the definition of $G$ it follows that both the inclusion $i:P\to T(P)$ and $F:P\mto T(P)$ are selectors of $G$, that is, for any $x\in P_1$ we have
$$
x\in G(x)\hspace{1cm}\mbox{ and }\hspace{1cm} F(x)\subset G(x).
$$
Moreover, all of the above maps are acyclic-valued (cf.\ again Proposition \ref{lem:D_properties} and Theorem~\ref{th:F-usc-acyclic}). Therefore, it follows from \cite[Proposition 32.13(i)]{Go06} that the identities
$
H^*(F)=H^*(G)=H^*(i)
$
are satisfied.
As a consequence we obtain the desired equality $I_P=\id_{H^*(P)}$, which completes the proof.
\qed\\

\section{Proof of Theorem~\ref{thm:ConleyF=ConleyV}}
\label{sec:ConleyF=ConleyV}

In order to prove Theorem~\ref{thm:ConleyF=ConleyV} we first construct an auxiliary pair $(Q_1,Q_2)$
and show that $H^*(P_1,P_2)\cong H^*(Q_1,Q_2)$. As a second step, we then construct a continuous surjection
$\psi:(Q_1,Q_2)\to (|\ccl\cS|,|\Exit\cS|)$ with contractible preimages and apply the Vietoris-Begle theorem
to complete the proof.

\subsection{The pair $(Q_1,Q_2)$.}
\begin{figure}[tb]
  \includegraphics[width=0.60\textwidth]{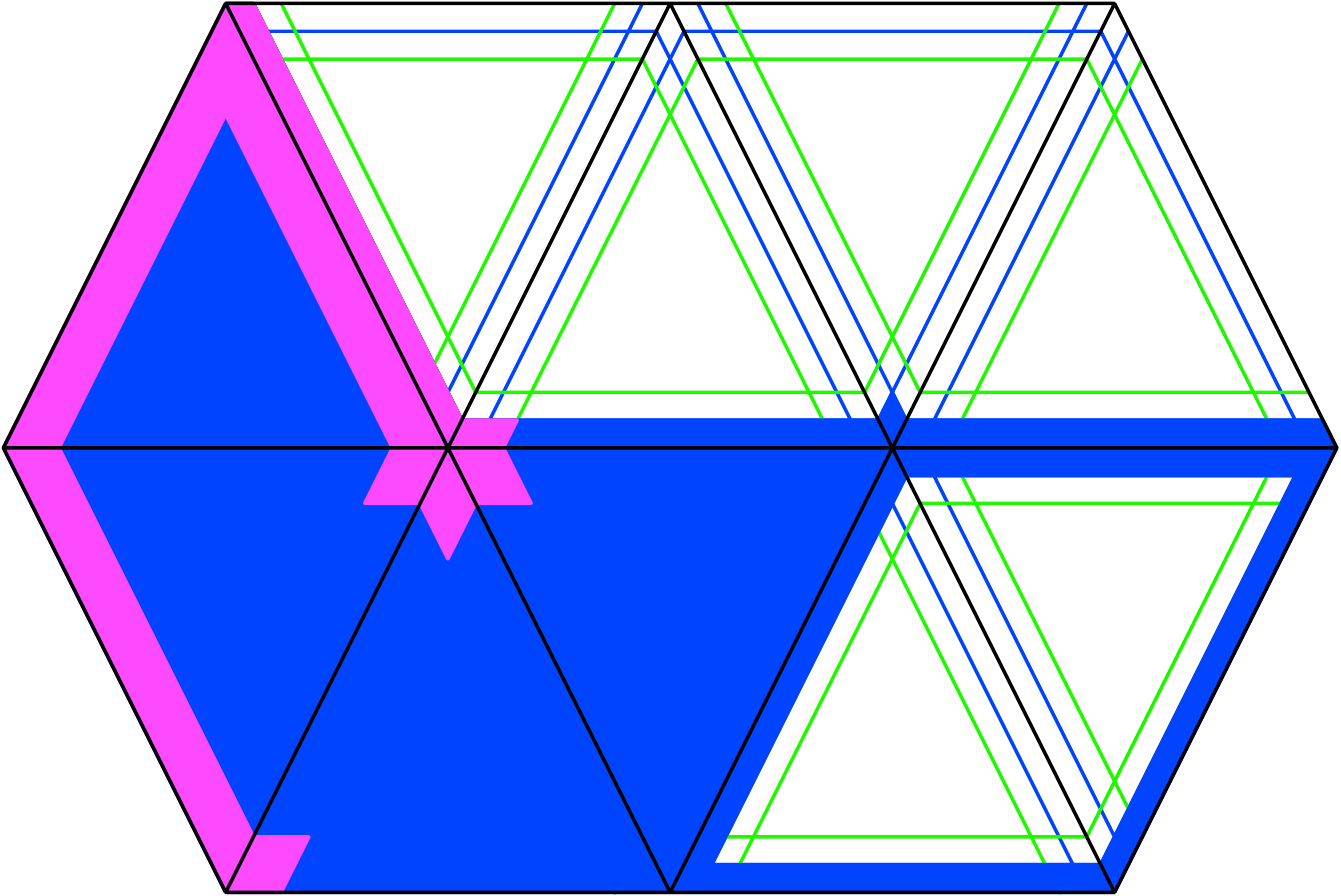}
  \caption{The pair~$Q = (Q_1,Q_2)$ associated with the isolating
           block~$N_\delta$ from Figure~\ref{fig:isoblock}
           and the weak index pair~$P = (P_1,P_2)$ from
           Figure~\ref{fig:isopairP}. The set~$Q_1$ is the union
           of the dark blue and magenta regions, while the subset
           $Q_2 \subset Q_1$ is only the magenta part.}
  \label{fig:isopairQ}
\end{figure}

Consider the pair $(Q_1,Q_2)$ consisting of the two sets
$$
\begin{array}{l}
Q_1:=N_\delta (\ccl \cS)\cap N_{\delta'} (\ccl \cS),\\[1ex]
Q_2:=N_\delta (\Exit \cS)\cap N_{\delta'} (\ccl \cS),
\end{array}
$$
where $N_\delta(\cA)$ is given by \eqref{eq:N-beta}.
Figure~\ref{fig:isopairQ} shows an example of such a pair for the isolated invariant set
$\cS$ presented in Figure~\ref{fig:isoblock}.
\begin{prop}\label{prop:P_iso_Q}
We have
\begin{equation}
\label{eq:P_iso_Q}
H^*(Q_1,Q_2) \cong H^*(P_1,P_2).
\end{equation}

\end{prop}
\proof We begin by verifying the two inclusions
\begin{equation}\label{eq:P2subsetQ2}
P_i\subset Q_i
\quad\mbox{ for }\quad
i=1,2.
\end{equation}
It is clear that $P_1\subset Q_1$, therefore we shall verify (\ref{eq:P2subsetQ2}) for $i=2$.

Let $x\in P_2$. Then $x\in N_\delta$ and by Lemma~\ref{lem:N-bd-cK}, there exists a simplex $\sigma\notin \cS$ such that $x\in\cl\csd{\sigma}$. On the other hand, $x\in P_1$ implies $x\in N_{\delta'}(\cS)$, so we can take $\tau\in\cS$ with $x\in\cl\csdp{\tau}$. For any vertex $v\in\sigma$ we have $t_v(x)\geq\delta>\delta'$, which shows that $\sigma\subset\tau\in\cS$. Consequently, $\sigma\in\ccl\cS$, which along with $\sigma\notin\cS$ implies
the inclusions $\sigma\in\Exit\cS$ and $x\in N_\delta(\Exit\cS)$. Observe now that we also have $x\in N_{\delta'}(\ccl\cS)$, according to $P_2\subset N_{\delta'}(\cS)\subset N_{\delta'}(\ccl\cS)$. Thus, $x\in Q_2$. The proof of (\ref{eq:P2subsetQ2}) is now complete.

Note that $P_1,P_2,Q_1,Q_2$ are compact and  $Q_2\subset Q_1$ and $P_2\subset P_1$.
Therefore, by the strong excision property of Alexander-Spanier cohomology, in order to
prove~\eqref{eq:P2subsetQ2}, it suffices to verify that
$
Q_1\setminus Q_2=P_1\setminus P_2.
$

For this, consider an $x\in Q_1\setminus Q_2$. Then $x\in N_\delta (\ccl\cS)$ and $x\notin
N_\delta (\Exit \cS)$.
Hence, there exists a $\sigma\in\ccl\cS\setminus\Exit\cS=\cS$ with $x\in\cl\csd{\sigma}$.
It follows that \begin{equation}\label{eq:x_in_Ndelta}
x\in N_\delta(\cS).
\end{equation}
We also have $x\in N_{\delta'}(\ccl\cS)$. In order to show that $x\in P_1$, we need to verify
the inclusion $x\in N_{\delta'}(\cS)$. Suppose to the contrary
that there is a $\tau\in\ccl\cS\setminus\cS=\Exit\cS$ with $x\in\cl\csdp{\tau}$. Then for each vertex $v$ of $\sigma$ we have $t_v(x)\geq \delta>\delta'$, which means that each vertex of $\sigma$ is a vertex of $\tau$.
In other words $\sigma\subset\tau$. However, $\tau\in\Exit\cS$ which, according to the closedness of $\Exit\cS$, implies $\sigma\in\Exit\cS$, a contradiction.
Therefore, $x\in N_{\delta'}(\cS)$, and together with~(\ref{eq:x_in_Ndelta}) this implies the inclusion $x\in P_1$.
Since $x\notin Q_2$, by (\ref{eq:P2subsetQ2}), we further have $x\notin P_2$. Consequently, both $x\in P_1\setminus P_2$
and $Q_1\setminus Q_2\subset P_1\setminus P_2$ are satisfied.

In order to prove the reverse inclusion let $x\in P_1\setminus P_2$ be arbitrary.
It is clear that then $x\in Q_1$. We need to show that $x\notin Q_2$. Suppose the contrary.
Then there exists a simplex $\sigma\in\Exit\cS$ such that $x\in\cl\csd{\sigma}$.
It follows that  $\sigma\in\cX^\delta(x)\setminus \cS$ and $\cX^\delta(x)\setminus \cS\neq\emptyset$.
Since $x\in P_1\subset N_\delta(\cS)$, we also have $\cX^\delta(x)\cap \cS\neq\emptyset$.
Therefore, by Lemma~\ref{lem:N-bd-cK}, we get $x\in\bd N_\delta(\cS)$.
Yet, we also have  $x\in N_{\delta'}(\cS)$ in view of $x\in P_1$. Consequently, $x\in   \bd N_\delta(\cS)\cap N_{\delta'}(\cS)=P_2\subset Q_2$, which is a contradiction.
\qed\\

\subsection{Auxiliary maps $\varphi_\sigma ^\lambda$.}
\begin{prop}\label{prop:intersection}
For any $\sigma,\tau\in\cX$ we have
$$
\begin{array}{lccl}
\cl\csl{\sigma}\cap\cl\csl{\tau}=\left\{x\in X\;\mid \right.&t_v(x)=\lambda&\mbox{ for }&v\in(\tau\setminus\sigma)\cup(\sigma\setminus\tau),\\[1ex]
&t_v(x)\geq\lambda&\mbox{ for }&v\in\tau\cap\sigma,\\[1ex]
&t_v(x)\leq\lambda&\mbox{ for }&v\notin\tau\cup\sigma\left.\right\}.
\end{array}
$$
In particular
$$
\cl\csl{\sigma}\cap\cl\csl{\tau}\subset\cl\csl{\sigma\cap\tau}.
$$
\end{prop}
\proof The proposition follows immediately from (\ref{eq:cl_sigma_l}).\qed\\
\ \ \

For $\lambda\in[0,1)$ let
\begin{equation}\label{phi_lambda}
\varphi _\lambda:[0,1]\ni t\;\longmapsto\;\left\{\begin{array}{cl}
\frac{t-\lambda}{1-\lambda}&\mbox{ for }t\geq\lambda\\[1ex]
0&\mbox{ for }t\leq\lambda
\end{array}
\right.\in[0,1].
\end{equation}
Given a simplex $\sigma$ in $\cX$ we define the map
$$
\varphi_\sigma ^\lambda:\cl\csl{\sigma}\ni x\longmapsto\varphi_\sigma ^\lambda(x) \in |\sigma|
$$
by
\begin{equation}
\label{phi_sigma}
\varphi_\sigma ^\lambda(x):=\sum_{v\in\cX_0}\frac{\varphi_\lambda (t_v(x))}{\sum_{w\in \cX_0 } \varphi_\lambda (t_w(x))}v\;.
\end{equation}
\begin{prop}\label{prop:lip}
The map $\varphi_\sigma^\lambda$ is
well-defined and continuous.
\end{prop}
\proof
Let $x\in\cl\csl{\sigma}$.
Then we have $t_v(x)\leq \lambda$ for $v\not\in\sigma$,
and consequently the identity
$\varphi_\lambda(t_v(x))=0$ holds
for $v\not\in\sigma$.
Hence, $\varphi_\sigma^\lambda(x)\in|\sigma|$, which means
that~$\varphi_\sigma ^\lambda$ is well-defined.
The continuity of $\varphi_\sigma ^\lambda(x)$
follows from the continuity of $\varphi_\lambda$ and the continuity of the barycentric coordinates.
\qed

For $\sigma\in\cX$ let $n_\sigma$ denote the number of vertices in $\sigma$. For $x\in\cl\csl{\sigma}$
set
$$
   r^\lambda_\sigma(x):=\sum_{w\not\in\sigma}t_w(x).
$$
\begin{lem}\label{lem:phisigma}
Let $\sigma\in\cX$ and $\lambda\in[0,1)$. For any $x\in \cl\csl{\sigma}$ and $v\in\cX_0$
we have
\begin{equation}\label{eq:phisigma}
t_v(\varphi_\sigma ^\lambda(x))=
\left\{\begin{array}{cl}\displaystyle{\frac{t_v(x)-\lambda}{1-\lambda n_\sigma-r^\lambda_\sigma(x)}} &
\mbox{\rm if }v\in\sigma\\[2ex]
0&\mbox{\rm otherwise}.
\end{array}\right.
\end{equation}
\end{lem}
\proof
It is clear from \eqref{phi_lambda} and \eqref{phi_sigma}
that (\ref{eq:phisigma}) is correct  for $v\notin\sigma$.
On the other hand,
if $v\in\sigma$ then $t_v(x)\geq\lambda$, hence, by (\ref{phi_sigma}) and (\ref{phi_lambda}) we have
\begin{equation}
\label{eq:phisigma-2}
t_v(x)=t_v(\varphi_\sigma ^\lambda(x))(1-\lambda)\sum_{w\in \cX_0 }\varphi_\lambda(t_w(x))+\lambda.
\end{equation}
Summing up the barycentric coordinates of~$x$ over all vertices in~$\cX_0$,
and taking into account the above equalities, which are valid for all vertices of $\sigma$,
we obtain
\begin{eqnarray*}
\sum_{v\in \cX_0 }t_v(x) & = & \sum_{v\in \sigma }t_v(x)+\sum_{v\not\in \sigma }t_v(x) \\[1.5ex]
 & = &
\sum_{v\in \sigma }t_v(\varphi_\sigma ^\lambda(x))(1-\lambda)\sum_{w\in \cX_0 }\varphi_\lambda(t_w(x))+\lambda n_\sigma+r^\lambda_\sigma(x).
\end{eqnarray*}
Since the barycentric coordinates sum to~$1$, we have $\sum_{v\in \cX_0 }t_v(x)=1$.
Moreover, since $\varphi_\sigma ^\lambda(x)\in|\sigma|$, we also have
$\sum_{v\in \sigma }t_v(\varphi_\sigma ^\lambda(x))=1$.
Therefore,  the above equality reduces to
$$
1=(1-\lambda)\sum_{w\in \cX_0 }\varphi_\lambda(t_w(x))+\lambda n_\sigma+r^\lambda_\sigma(x).
$$
Consequently,
$$
\sum_{w\in \cX_0 }\varphi_\lambda(t_w(x))=\frac{1-\lambda n_\sigma-r^\lambda_\sigma(x)}{1-\lambda}.
$$
Replacing the sum  $\sum_{w\in \cX_0 }\varphi_\lambda(t_w(x))$ in \eqref{eq:phisigma-2} by the right-hand
side of this equation
 and calculating $t_v(\varphi_\sigma ^\lambda(x))$
we obtain \eqref{eq:phisigma} for $v\in\sigma$.
This completes the proof.
\qed
\begin{prop}\label{prop:phi_1}
For any simplex $\sigma\in\cX$ and $\lambda\in[0,1)$ we have
$$
\varphi_\sigma ^\lambda(\sigma\cap\cl\csl{\sigma})=\sigma.
$$
\end{prop}
\proof It is clear that $
\varphi_\sigma ^\lambda(\sigma\cap\cl\csl{\sigma})\subset\sigma,
$
therefore we verify the opposite inclusion. Take an arbitrary $y\in\sigma$ and define
$$
x:=\sum_{v\in\sigma}\left(t_v(y)(1-\lambda n_\sigma)+\lambda\right)v
$$
It is easy to check that the above formula correctly defines a point $x\in \sigma$ via its barycentric coordinates. Moreover, we have $x\in\sigma\cap\cl\csl{\sigma}$, as $t_v(x)=0$ for $v\notin\sigma$ and $t_v(x)\geq\lambda$ for $v\in\sigma$.
An easy calculation, with the use of Lemma \ref{lem:phisigma}, finally shows that $\phi_\sigma ^\lambda(x)=y$.
\qed

\medskip

The following proposition is an immediate consequence of Proposition \ref{prop:intersection}.
\begin{prop}\label{prop:phi_intersect}
For any simplices $\sigma$ and $\tau$ in $\cX$ and arbitrary $\lambda\in[0,1)$
the maps~$\varphi_\sigma^\lambda$ and~$\varphi_\tau^\lambda$ coincide on
$\cl\csl{\sigma}\cap \cl\csl{\tau}$
\qed
\end{prop}

\subsection{Mapping $\psi$.}
In view of~\eqref{eq:x} and Proposition~\ref{prop:phi_intersect}, we have a well-defined continuous surjection $\varphi:|\cX|\to |\cX|$ given by
$$
\varphi(x):=\varphi_\sigma^\delta (x)
\quad\mbox{ where $\sigma\in\cX$ is such that }\quad
x\in \cl\csl{\sigma}.
$$
Let $\psi:=\varphi_{|Q_1}:Q_1\to X$ denote the restriction of $\varphi$ to $Q_1$.

\begin{prop}\label{prop:fibers}
For each $y\in|\ccl\cS|$ the fiber $\psi^{-1}(y)$
is non-empty and contractible.
\end{prop}

\proof
Let $y\in |\ccl\cS|=\ogr{\ccl\cS}$ be arbitrary and let the simplex $\sigma\in\ccl\cS$ be such that $y\in\cell\sigma$. Furthermore, define the set
$$
\begin{array}{ll}X_\sigma:=\left\{\right.x\in X\;\mid
&t_v(x)\leq\delta\mbox{ if }v\notin\sigma\mbox{ and }\\[1ex]
&
t_v(x)=t_v(y)(1-\delta n_\sigma-r^\lambda_\sigma(x))+\delta
\mbox{ for }v\in\sigma\left.\right\}.
\end{array}
$$
We first verify that the fiber of $y$ under $\varphi _\sigma ^\delta$ is given by
\begin{equation}\label{eq:counterimage_phi_sigma}
(\varphi _\sigma ^\delta) ^{-1}(y)=X_\sigma.
\end{equation}
For this, fix an  $x\in X_\sigma$, and
recall that $\delta$ satisfies (\ref{def:deltagammaepsbound}), in particular $\delta<1/(d+1)$.
Therefore, for $v\in\sigma$  we deduce from \eqref{eq:phisigma} the inequality
$$
t_v(x)=t_v(y)(1-\delta n_\sigma-r^\lambda_\sigma(x))+\delta
\geq t_v(y)(1-\delta n_\sigma-\delta(d+1-n_\sigma)+\delta\geq \delta.
$$
This, together with the obvious inequality $t_v(x)\leq\delta$ for $v\notin\sigma$, shows that
\begin{equation}\label{eq:phi-1_incl_clsigma}
X_\sigma\subset\cl\csd{\sigma}=\dom \varphi_\sigma^\delta.
\end{equation}
Moreover, a straightforward calculation implies that for every point
$x\in X_\sigma$ the identity $\varphi_\sigma^\delta(x)=y$ holds.
This shows that $X_\sigma\subset (\varphi _\sigma ^\delta)^{-1}(y)$. Since the converse inclusion is straightforward, the proof of (\ref{eq:counterimage_phi_sigma}) is finished.
Now let
\begin{equation}\label{eq:Ydef}
\bar{X}_\sigma:=(\varphi _\sigma ^\delta)^{-1}(y)\cap  N_{\delta'}(\ccl\cS)=X_\sigma\cap N_{\delta'}(\ccl\cS).
\end{equation}
We claim that
\begin{equation}\label{eq:psi-1}
\psi^{-1}(y)=\bar{X}_\sigma.
\end{equation}
Note that
\[
  \bar{X}_\sigma\subset \cl\csd{\sigma}\cap N_{\delta'}(\ccl\cS)\subset  N_{\delta}\cap  N_{\delta'}=Q_1=\dom \psi.
\]
Since $\bar{X}_\sigma\subset(\varphi_\sigma ^\delta) ^{-1}(y)$, one obtains for $w\in \bar{X}_\sigma$ the identity
$y=\varphi_\sigma ^\delta(w)=\psi(w)$. Therefore, $\bar{X}_\sigma\subset \psi^{-1}(y)$.
For the proof of the reverse inclusion it suffices to verify that the condition $x\notin\cl\csd{\sigma}$ implies $\psi(x)\neq y$.
Suppose to the contrary that $x\not\in\cl\csd{\sigma}$ and $\psi(x)=y$, and
consider a simplex $\tau$ such that $x\in \cl\csd{\tau}$.
Then  $\psi(x)=\varphi_\tau ^\delta (x)$.
Directly from the definition of $\varphi_\tau ^\delta$ we infer that $\psi(x)\in\tau$. However $\psi(x)=y\in\cell\sigma$, which means that $\sigma$ is a face of $\tau$. Then, taking into account the inclusion $x\in\cl\csd{\tau}\setminus \cl\csd{\sigma}$, we can find a vertex $v\in\tau\setminus\sigma$ such that $t_v(x)>\delta$. Consequently, by (\ref{eq:phisigma}), we have $t_v(y)=t_v(\psi(x))=t_v(\varphi_\tau ^\delta (x))>0$, which contradicts $y\in\cell\sigma$, and completes the proof of (\ref{eq:psi-1}).

We still need to show that $\bar{X}_\sigma=\psi^{-1}(y)$ is contractible. To this end, we define the map
$
h:\bar{X}_\sigma\times[0,1]\to \bar{X}_\sigma
$
by
\[
   h(x,s):=\sum_{v\in\cX_0} t_{v,x,s}v,
\]
where
$$
t_{v,x,s}:=\left\{\begin{array}{ll}
(1- s)t_v(x)&\mbox{ if }v\notin\sigma,\\[1ex]
t_v(y)(1-\delta n_\sigma-(1- s)r^\lambda_\sigma(x))+\delta&\mbox{ if }v\in\sigma.
\end{array}\right.
$$
We will show that $h$  is a well-defined homotopy between the identity on $\bar{X}_\sigma$ and a constant map on $\bar{X}_\sigma$.

To begin with, we verify that  for any point $x\in \bar{X}_\sigma$ and arbitrary $s\in[0,1]$ we have
$h(x, s)\in \bar{X}_\sigma$. The verification that the inclusion $h(x, s)\in\cl\csd{\sigma}$ holds, as well as $\varphi_{\sigma} ^{\delta}(h(x, s))=y$, which in turn shows that $h(x, s)\in(\varphi_\sigma ^\delta)^{-1}(y)$, is tedious but straightforward.
We still need to verify that $h(x, s)\in N_{\delta'}(\ccl\cS)$.
For this, consider a simplex $\tau\in\ccl\cS$ such that $x\in\cl\csd{\sigma}\cap\cl\csdp{\tau}$.
Since for any $v\in\sigma$ we have $t_v(x)\geq\delta>\delta'$, we deduce that $\sigma\subset\tau$.
Let
\[
   \eta:=\setof{v\in\cX_0\mid t_v(h(x,s))>\delta'}.
\]
Then we claim that the inclusions $\sigma\subset\eta\subset\tau$ hold.
Indeed, if $v\notin\tau$ then $\sigma\subset\tau$ implies the inequalities $t_v(h(x, s))=(1- s)t_v(x)\leq t_v(x)\leq\delta'$.
Furthermore, if one has $v\in\sigma\subset\tau$, then $t_v(h(x, s))\geq t_v(x)\geq \delta>\delta'$.
Therefore, $\eta$ is a simplex and it satisfies $h(x, s)\in\cl\csdp{\eta}$.
Since $\sigma\in\ccl\cS$ as  well as $\tau\in\ccl\cS$, the closedness of $\ccl\cS$ implies
that $\eta\in\ccl\cS$.
Consequently, $h(x, s)\in N_{\delta'}(\ccl\cS)$, and this proves that the map~$h$ is well-defined.

The continuity of $h$~follows from the continuity of the barycentric coordinates. Verification that $h(\cdot,0)=\id_{\bar{X}_\sigma}$ as well as that $h(\cdot,1)$ is constant on $\bar{X}_\sigma$ is straightforward. This completes the proof.
\qed
\begin{prop}\label{prop:psi_surj}
We have
\begin{itemize}
\item[(i)] $\psi(Q_1)=|\ccl\cS|$,
\item[(ii)] $\psi(Q_2)=|\Exit\cS|$,
\item[(iii)] $\psi ^{-1}(|\Exit\cS|)=Q_2$.
\end{itemize}
In particular, we can consider $\psi$ as a map of pairs
$$\psi:(Q_1,Q_2)\to (|\ccl\cS|,|\Exit\cS|).$$
 \end{prop}
\proof For the proof of (i), fix an arbitrary point $x\in Q_1$. Then there exists a simplex $\sigma\in\ccl\cS$
with $x\in\cl\csd{\sigma}$.
Thus, $\psi(x)=\varphi(x)=\varphi_\sigma ^\delta(x)\in\sigma\subset|\ccl\cS|$ and  $\psi (x)\in|\ccl\cS|$.
This implies that the inclusion $\psi(Q_1)\subset|\ccl\cS|$.
The reverse inclusion is a consequence of Proposition \ref{prop:phi_1},
because for any simplex $\sigma\in\ccl\cS$ we have
$\sigma\cap\cl\csd{\sigma}\subset\cl\csd{\sigma}\cap\cl\csdp{\sigma}\subset Q_1$.
The proof of (ii) is  analogous to the proof of (i).

In order to prove the remaining statement (iii), first observe that we have
the inclusion $\psi ^{-1}(|\Exit\cS|)\subset Q_2$. Indeed, given a
$y\in|\Exit\cS|$ there exists a $\sigma\in\Exit\cS$ such that
$y\in\cell\sigma$, and by (\ref{eq:Ydef}), (\ref{eq:psi-1}),
and~(\ref{eq:phi-1_incl_clsigma}), we have
$$
\psi ^{-1}(y)\subset\cl\csd{\sigma}\cap N_{\delta'}(\ccl\cS)\subset N_{\delta}(\Exit\cS)\cap N_{\delta'}(\ccl\cS)=Q_2,
$$
which implies $\psi ^{-1}(|\Exit\cS|)\subset Q_2$. This, together with (ii), implies (iii).
The last statement is a direct consequence of (i) and (ii).\qed
\begin{prop}\label{prop:Q_iso_ClS}
We have
$$
H^*(Q_1,Q_2) \cong H^*(|\ccl\cS|,|\Exit \cS|).
$$
\end{prop}
\begin{figure}[tb]
\begin{minipage}[t]{0.47\linewidth}
\centering
\includegraphics[width=\textwidth]{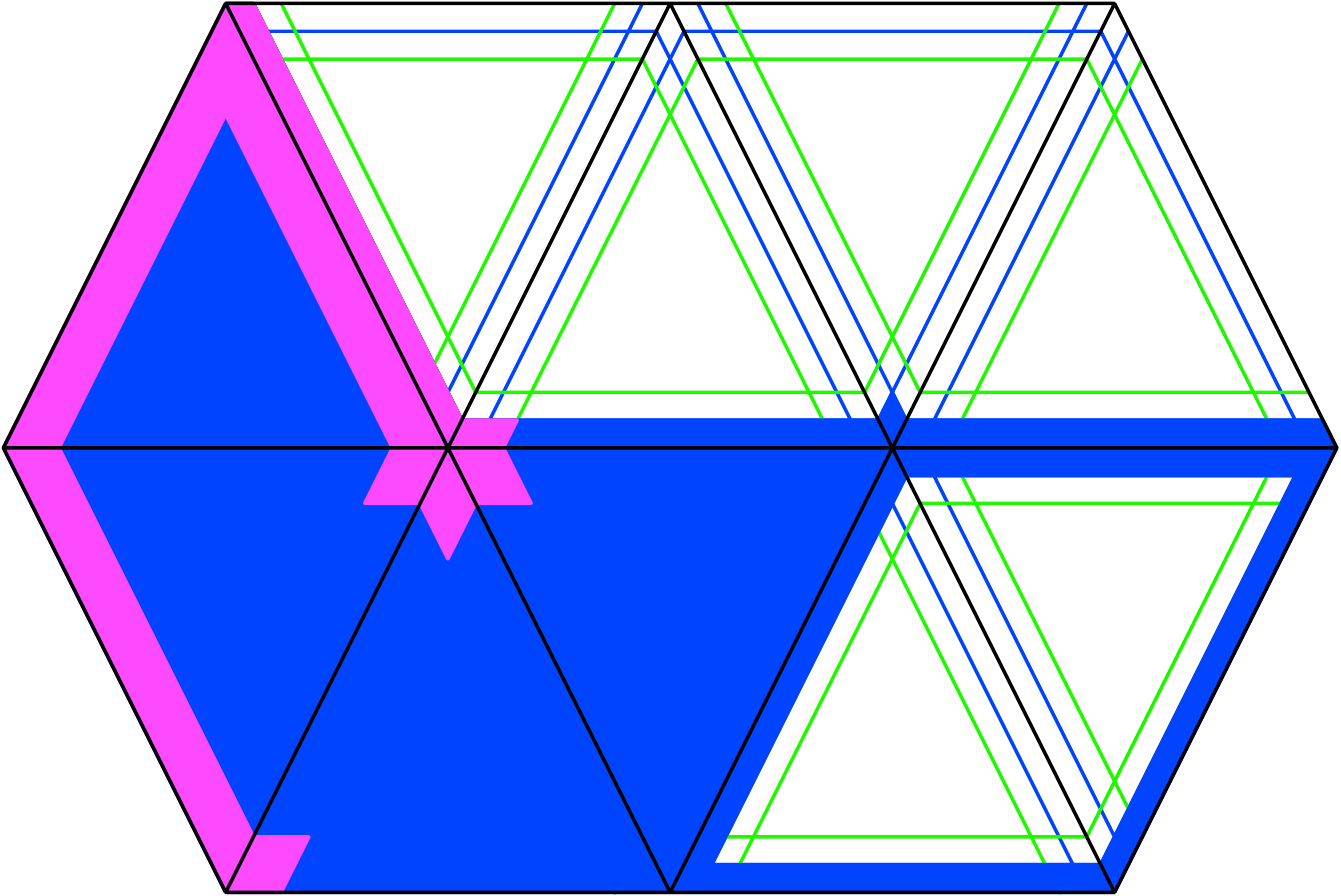}
\subcaption{The pair $Q=(Q_1,Q_2)$. $Q_1$ is the union of the dark blue
and magenta regions, while~$Q_2$ is only the magenta region.}
\end{minipage}
\centering
\hspace{5mm}\begin{minipage}[t]{0.47\linewidth}
\includegraphics[width=\textwidth]{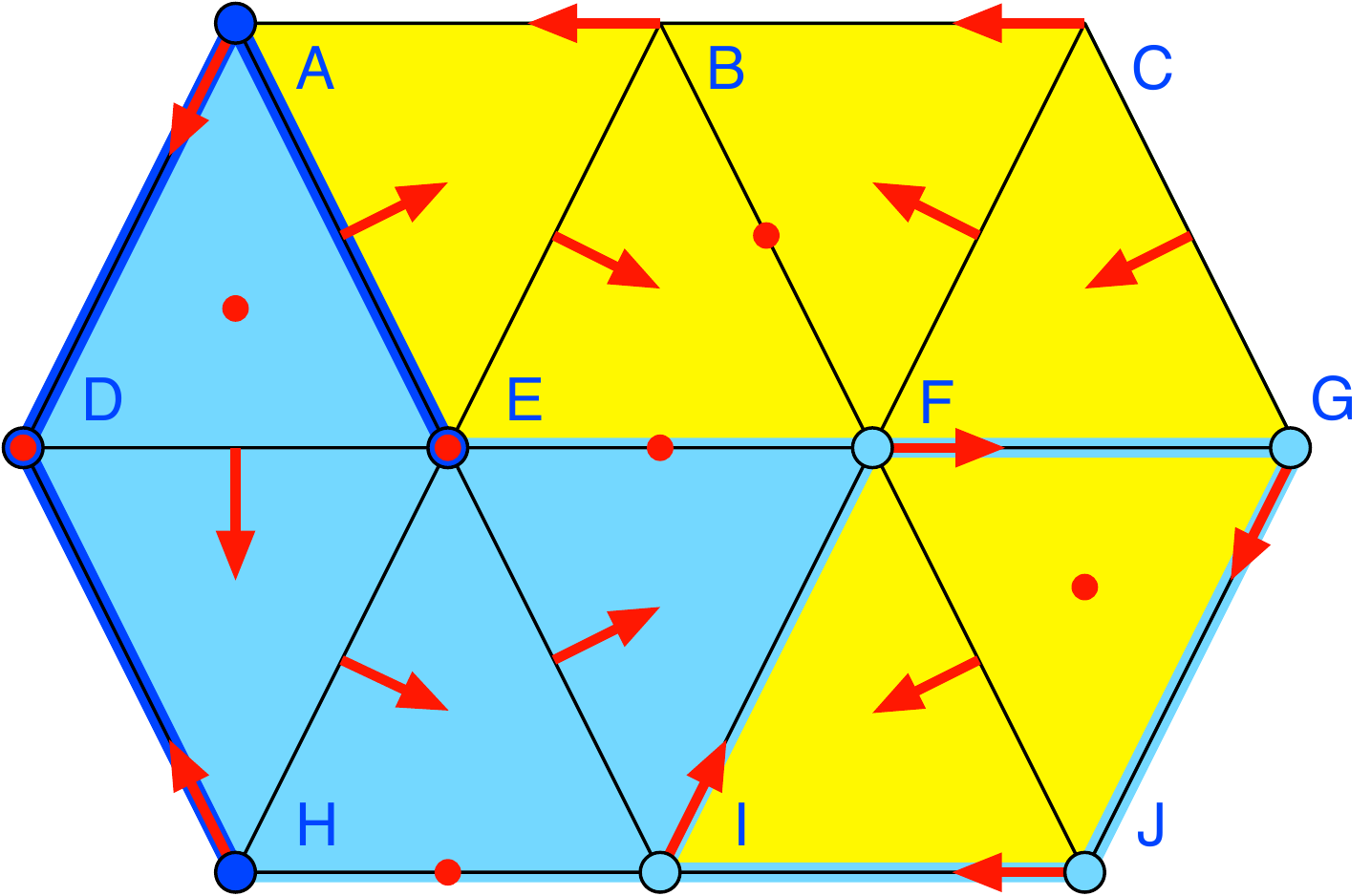}
\subcaption{The pair $(|\ccl\cS|,|\Exit\cS|)$. $|\ccl\cS|$ is the union of
the light blue and dark blue regions, while~$|\Exit\cS|$ is only the dark
blue set.}
\end{minipage}
\caption{The pairs $Q=(Q_1,Q_2)$ and $(|\ccl\cS|,|\Exit\cS|)$.}
\label{fig:pairs_Q_and_clsexit}
\end{figure}

\proof
By Proposition \ref{prop:psi_surj} the mapping $\psi:(Q_1,Q_2)\to (|\ccl\cS|,|\Exit\cS|)$ is a continuous surjection  with $\psi ^{-1}(|\Exit\cS|)=Q_2$. By Proposition \ref{prop:fibers}, $\psi$ has contractible, and hence acyclic fibers.
Moreover, $\psi$ is proper, that is, the counterimages of compact sets under~$\psi$ are compact.
Therefore, the map~$\psi$ is a Vietoris map. By the Vietoris-Begle mapping theorem for the pair of spaces we conclude that
$$
\psi ^*:H^*(|\ccl\cS|,|\Exit\cS|)\to
H^*(Q_1,Q_2)$$
is an isomorphism, which completes the proof.\qed
\\
\ \ \

Figure \ref{fig:pairs_Q_and_clsexit} shows an example of the pairs
$(Q_1,Q_2)$ and $(|\ccl\cS|,|\Exit \cS|)$ in Proposition \ref{prop:Q_iso_ClS}.\\
\ \ \

{\bf Proof of Theorem~\ref{thm:ConleyF=ConleyV}:} Theorem~\ref{thm:ConleyF=ConleyV} is an immediate consequence of
Proposition~\ref{prop:P_iso_Q}, Proposition~\ref{prop:Q_iso_ClS}, and Theorem~\ref{thm:F-has-flowlike-Conley-index}.
\qed

\section{Proof of Theorem~\ref{thm:Mcomb=MF}}
\label{sec:Mcomb=MF}

In order to prove Theorem~\ref{thm:Mcomb=MF} we first establish a few auxiliary lemmas.
Then we recall some results concerning the correspondence of solutions
for~$\Pi_\cV$ and~$F$. We then use this correspondence to prove an auxiliary theorem
and finally present the proof of Theorem~\ref{thm:Mcomb=MF}.

\subsection{Auxiliary lemmas.}
First observe that Theorem~\ref{th:N-iso-block} applies to the set $N_\beta(\cS)$ given by \eqref{eq:N-beta}
for any $\beta$ which satisfies $0<\beta< 1 / (d+1)$.

\begin{lem}\label{lem:Ne_cap_supS_closed}
We have
$$
N_\epsilon\cap \ogr{\cS}=N_\epsilon\cap |\ccl\cS|.
$$
In particular, $
N_\epsilon\cap \ogr{\cS}$ is closed.
\end{lem}
\proof
Clearly $N_\epsilon\cap \ogr{\cS}\subset N_\epsilon\cap |\ccl\cS|$.
To prove the opposite inclusion, assume to the contrary that
there exists an $x\in N_\epsilon\cap |\ccl\cS|$ and $x\not\in N_\epsilon\cap \ogr{\cS}$.
Then, by Proposition~\ref{prop:ogr}(ii), $x\in N_\epsilon\cap |\Exit\cS|$.
Consider simplices $\sigma\in\Exit\cS$ and $\tau\in\cS$ such that $x\in\cell\sigma$ and $x\in\cl\cse{\tau}$.
Since for any vertex $v\in \tau$ we have $t_v(x)\geq\varepsilon>0$, the inclusion $v\in\sigma$
has to hold. Hence, $\tau\subset\sigma$. Therefore, by the closedness of $\Exit\cS$ we get
$\tau\in\Exit\cS$, a contradiction.
\qed
\begin{lem}\label{lem:sigma_max_in_S}
For any $x\in N_\epsilon\cap \ogr{\cS}$ we have
$
\sigma^\epsilon_{max}(x)\in\cS.
$
\end{lem}
\proof
Fix a point $x\in N_\epsilon\cap \ogr{\cS}$. Then there exist simplices $\tau,\sigma\in\cS$ such that $x\in\cl\cse{\tau}$ and $x\in\cell{\sigma}$. Clearly, one has $\sigma=\sigma^0_{min}(x)$.
By Corollary~\ref{cor:equiv-cell-closure} and Lemma~\ref{lem:epsilon-delta-include} we then obtain
$
\tau\subset\sigma^\epsilon_{max}(x)\subset\sigma.
$
Therefore, the closedness of $\Exit\cS$ implies that $\sigma^\epsilon_{max}(x)\in\cS$.
\qed
\begin{lem}\label{lem:NepsilonsubsetNdelta}
We have
$$
N_\epsilon\cap \ogr{\cS}\subset N_\delta\cap\ogr{\cS}.
$$
\end{lem}
\proof Fix a point $x\in N_\epsilon\cap \ogr{\cS}$.
Then we have $\sigma^0_{min}(x)\in\cS$. By Lemma \ref{lem:sigma_max_in_S} we further obtain $\sigma^\epsilon_{max}(x)\in \cS$.
In addition, Lemma~\ref{lem:epsilon-delta-include} immediately implies the inclusions
$
\sigma^\epsilon_{max}(x)\subset \sigma^\delta_{min}(x)\subset\sigma^\delta_{max}(x)\subset \sigma^0_{min}(x).
$
Now the closedness of $\Exit\cS$ yields $\sigma^\delta_{max}(x)\in\cS$,
and consequently, $x\in N_\delta$, which completes the proof.\qed

\subsection{Solution correspondence.}
In the sequel we need two results on the correspondence of solutions of
the combinatorial flow $\Pi_\cV$ and the associated multivalued dynamical system $F$.
We recall them from \cite{KaMrWa16}. We begin with a definition.
\begin{defn} \label{defn:redextsoln} (see \cite[Definition 5.2]{KaMrWa16})
\begin{itemize}
\item[(a)] Let~$\rho : \ZZ \to \cX$ denote a full solution
of the combinatorial flow~$\Pi_\cV$. Then the {\em reduced
solution\/}~$\rho^* : \ZZ \to \cX$ is obtained
from~$\rho$ by removing~$\rho(k+1)$ whenever~$\rho(k+1)$
is the target of an arrow of~$\cV$ whose source is~$\rho(k)$.
\item[(b)] Conversely, let~$\rho^* : \ZZ \to \cX$ denote
an arbitrary sequence of simplices in~$\cX$. Then its {\em
arrowhead extension\/}~$\rho : \ZZ \to \cX$ is defined as
follows. If~$\rho^*(k) \in \dom\cV \setminus \Fix\cV$ and
if~$\rho^*(k+1) \neq \rho^*(k)^+$, then we insert~$\rho^*(k)^+$
between~$\rho^*(k)$ and~$\rho^*(k+1)$. In other words, the
arrowhead extension~$\rho$ is obtained from~$\rho^*$ by
inserting missing targets of arrows.
\end{itemize}
\end{defn}
\begin{thm} \label{thm:orbitequiv1}(see \cite[Theorem 5.3]{KaMrWa16})
Let
\begin{displaymath}
  X^\epsilon := \bigcup_{\sigma \in \cX} \cse{\sigma}
  \subset X
\end{displaymath}
denote the union of all open $\epsilon$-cells of~$X$.
Then the following hold.
\begin{itemize}
\item[(a)] Let~$\rho : \ZZ \to \cX$ denote a full solution
of the combinatorial flow~$\Pi_\cV$. Furthermore, let~$\rho^* :
\ZZ \to \cX$ denote the reduced solution as in
Definition~\ref{defn:redextsoln}(a). Then there is a
function~$\phi: \ZZ \to X^\epsilon$ such that for
$k \in \ZZ$ we have
\begin{displaymath}
  \phi(k+1) \in F\left( \phi(k) \right)
  \quad\mbox{ and }\quad
  \phi(k) \in \cse{\rho^*(k)} \; .
\end{displaymath}
In other words, $\phi$ is an orbit of~$F$ which follows
the dynamics of the combinatorial simplicial solution~$\rho$
after removing arrowheads.
\item[(b)] Conversely, let~$\phi : \ZZ \to X^\epsilon$
denote a full solution of~$F$ which is completely contained
in~$X^\epsilon$. Let $\rho^*(k) = \sigma_{max}^\epsilon(\phi(k))$
for $k \in \ZZ$, and let~$\rho : \ZZ \to \cX$ denote the
arrowhead extension of~$\rho$ as in Definition~\ref{defn:redextsoln}(b).
Then~$\rho$ is a solution of the combinatorial flow~$\Pi_\cV$.
\end{itemize}
\qed
\end{thm}

\begin{lem} \label{lem:F-sigma+} (see \cite[Lemma 4.9]{KaMrWa16})
For all simplices $\sigma \in \cX$ and all points $x \in X$ we have
$  F_\sigma(x) \subset \sigma^+$.
\qed
\end{lem}

\begin{thm} \label{thm:orbitequiv2} (see \cite[Theorem 5.4]{KaMrWa16})
Let~$\phi : \ZZ \to X$ denote an arbitrary full solution
of the multivalued map~$F$ and let
\begin{equation}
\label{eq:orbitequiv2-1}
\rho^*(k) = \sigma_{max}^\epsilon(\phi(k))
\end{equation}
for $k \in \ZZ$.
Extend this sequence of simplices in the following way:
\begin{itemize}
\item[(1)] For all $k \in \ZZ$ with $\phi(k)
\not\in |\rho^*(k-1)^+|$, we choose a face $\tau \subset
\rho^*(k-1)$ such that $\phi(k) \in |\ccl\tau^+
\setminus \{ \tau \}|$, and then insert~$\tau$
between~$\rho^*(k-1)$ and~$\rho^*(k)$.
\item[(2)] Let~$\rho : \ZZ \to \cX$ denote the
arrowhead extension of the sequence created in~{\em (1)},
according to Definition~\ref{defn:redextsoln}(b).
\end{itemize}
Then the so-obtained simplex sequence~$\rho : \ZZ \to \cX$
is a solution of the combinatorial flow~$\Pi_\cV$.
\qed
\end{thm}

\subsection{Invariance.}

\begin{lem}\label{lem:Nepsilon-negativelyinvariant}
The set $N_\epsilon\cap \ogr{\cS}$ is negatively invariant with respect to $F$, that is
$$
\Inv_F ^-(N_\epsilon\cap \ogr{\cS})=N_\epsilon\cap \ogr{\cS}.
$$
\end{lem}
\proof
Obviously, it suffices to prove that for every $y\in N_\epsilon\cap\ogr{\cS}$
there exists an $x\in N_\epsilon\cap \ogr{\cS}$ such that $y\in F(x)$.
To verify this, fix a $y\in N_\epsilon\cap \ogr{\cS}$.
Let $\sigma\in\cS$ be such that $y\in \cell\sigma$.
We will consider several cases concerning the simplex~$\sigma$.
First assume that  $\sigma\in \Fix \cV$, that is, $\sigma=\sigma^-=\sigma^+$. Take any $x\in\cse{\sigma}\cap\sigma\subset N_\epsilon\cap\ogr{\cS}$.
Since  $\sigma^\epsilon_{max}(x)=\sigma^-=\sigma^+$, the definition of $F$ (see \eqref{eq:defF}) shows that $F_\sigma (x)=\sigma$. Hence, $y\in F_\sigma(x)\subset F(x)$.

Now assume that $\sigma^-\neq\sigma^+$.
Note that if $\sigma=\sigma^+$ then we have $\sigma=A_\sigma\cup B_\sigma$, and if
$\sigma=\sigma^-$, then $\sigma\subset A_\sigma$.
Hence, either $y\in A_\sigma$
or $\sigma =\sigma^+$ and $y\in B_\sigma$.
In the latter case we may take
any point $x\in\cse{\sigma}\cap\sigma\subset N_\epsilon\cap\ogr{\cS}$,
because in that case one has
$\sigma=\sigma^\epsilon_{max}(x)=\sigma^\epsilon_{max}(x)^+\neq \sigma^\epsilon_{max}(x)^-$,
which immediately yields the inclusion $y\in B_\sigma=F_\sigma(x)\subset F(x)$.

It remains to consider the case $y\in A_\sigma$. Since $\cS$ is invariant with respect to $\Pi_\cV$ and $\sigma^+=\cV(\sigma^-)$, there exists a trajectory $\rho$ of $\cV$ in $\cS$ which contains $\sigma^-$ and $\sigma^+$ as consecutive simplices.
Let $\tau$ denote the simplex in this solution which precedes the tail $\sigma^-\in\dom\cV$.
Then $\tau\in\cS$ and, according to the definition of the multivalued flow $\Pi_\cV$,
we have $\sigma^-\subsetneq\tau\neq\sigma^+$.
Now let $k$ denote the number of vertices in $\tau\setminus\sigma$
and let $x\in X$ be the point with the barycentric coordinates given by
\[
    t_v(x):=\begin{cases}
              \epsilon & \text{if $v\in\tau\setminus\sigma$, }\\
              1-k\epsilon & \text{if $v\in\sigma$, }\\
              0 & \text{otherwise.}
            \end{cases}
\]
Then we have both $x\in\cl\cse{\tau}\cap\cl\cse{\sigma^-}\cap\tau$
and $\sigma^\epsilon_{max}(x)=\tau$, and this in turn implies $\sigma^\epsilon_{max}(x)^+\neq\sigma^-$ and $\sigma^\epsilon_{max}(x)^-\neq\sigma^-$. Therefore, $F_{\sigma^-}(x)=A_{\sigma^-}=A_\sigma$, which shows that $y\in F_{\sigma^-}(x)\subset F(x)$.
\qed

\begin{lem}\label{lem:Nepsilon-positivelyinvariant}
The set $N_\epsilon\cap \ogr{\cS}$ is positively invariant with respect to $F$, that is
$$
\Inv_F ^+(N_\epsilon\cap \ogr{\cS})=N_\epsilon\cap \ogr{\cS}.
$$
\end{lem}
\proof
For the proof it is enough to justify that for any point $x\in N_\epsilon\cap \ogr{\cS}$ we have $F(x)\cap N_\epsilon\cap \ogr{\cS}\neq\emptyset$. Let $x\in N_\epsilon\cap \ogr{\cS}$ be fixed and let $\sigma:=\sigma^\epsilon_{max}(x)$. By Lemma \ref{lem:sigma_max_in_S} we have $\sigma\in \cS$.  Then $x\in\cl\cse{\sigma}$. Since $F(\cse{\sigma})\subset F(\cl\cse{\sigma})$
and $F$ is strongly upper semicontinuous by Theorem~\ref{th:F-usc-acyclic},
without loss of generality we may assume that $x\in\cse{\sigma}$.
The set $\cS$ is invariant with respect to the combinatorial flow~$\Pi_\cV$.
Hence, there exists a solution $\rho:\ZZ\to \cX$ of $\Pi_\cV$, which is contained in $\cS$ and passes through~$\sigma$.
Furthermore, let $\rho^*:\ZZ\to\cX$ denote the reduced solution as defined in Definition~\ref{defn:redextsoln}(a).
There are two possible complementary cases: $\sigma\in\im\rho^*$ or $\sigma\notin\im\rho^*$.

In the first case there exists a $k\in\ZZ$ with $\sigma=\rho^*(k)$. Consider $\phi:\ZZ\to \cX^\epsilon$, which is a corresponding solution with respect to $F$ as constructed in Theorem~\ref{thm:orbitequiv1}(a). Then $\phi (k)\in\cse{\sigma}$ and $\phi (k+1)\in F(\phi(k))\cap\cse{\rho^*(k+1)}$. Since the map~$F$ is constant on open $\epsilon$-cells, we further obtain
\begin{equation}\label{eq:4.4e1}
\phi (k+1)\in F(x)\cap\cse{\rho^*(k+1)}.
\end{equation}
Due to $x\in\cse{\sigma}$, by Lemma~\ref{lem:F-sigma+} we have $F(x)=F_\sigma(x)\subset\sigma ^+\in\cS$, where the last inclusion follows from Proposition~\ref{prop:exSclosed+-}, as $\cS$ is an isolated invariant set. This, along with (\ref{eq:4.4e1}), completes the proof in the case where $\sigma\in\im\rho^*$.

Finally, we consider the case $\rho(k):=\sigma\notin\im\rho^*$, which immediately gives rise to the inclusion $\sigma\in\im\cV\setminus\Fix\cV$. In this case, the identity $\Pi_\cV (\sigma)=\Bd\sigma\setminus\{\sigma^-\}$
implies $\rho (k+1)\in \Bd\sigma\setminus\{\sigma^-\}$.
However, we also have $F(x)=F_\sigma (x)=B_\sigma$, according to the fact that $\sigma=\sigma^\epsilon _{max}(x) ^+\neq\sigma^\epsilon _{max}(x)^-$. This readily furnishes the inclusions $|\rho(k+1)|\subset|\Bd\sigma\setminus\{\sigma^-\}|\subset B_\sigma=F(x)$.
In particular, the barycenter of $\rho(k+1)$ belongs to $|\rho(k+1)|\cap\cse{\rho(k+1)}\cap F(x)\subset F(x)\cap N_\epsilon\cap\ogr{\cS}$. This completes the proof.
\qed
\medskip

As a straightforward consequence of Lemma \ref{lem:Nepsilon-negativelyinvariant} and Lemma \ref{lem:Nepsilon-positivelyinvariant} we obtain the following corollary.
\begin{cor}\label{cor:Nepsilon-invariant}
The set $N_\epsilon\cap \ogr{\cS}$ is invariant with respect to $F$, that is, we have
$$
\Inv_F (N_\epsilon\cap \ogr{\cS})=N_\epsilon\cap \ogr{\cS}.
$$
\qed
\end{cor}

\subsection{An auxiliary theorem and lemma.}
The following characterization of the set $\isoS=\Inv N_\delta$
is needed in the proof of Theorem~\ref{thm:Mcomb=MF}.
\begin{thm}\label{thm:invNd=invNe}
We have
$$
\Inv N_\delta=N_\epsilon\cap \ogr{\cS}.
$$
\end{thm}
\proof
According to Lemma \ref{lem:NepsilonsubsetNdelta} and Corollary~\ref{cor:Nepsilon-invariant}
we immediately obtain the inclusion $N_\epsilon\cap \ogr{\cS}\subset \Inv_FN_\delta$.
Therefore, it suffices to verify the opposite inclusion.

To accomplish this, take an $x\in N_\delta\setminus(N_\epsilon\cap \ogr{\cS})$.
If $x\in N_\delta\setminus N_\epsilon$ then $x\in\cl\csd{\tau}$ for some simplex $\tau\in\cS$, and $\sigma^\epsilon_{max}(x)\notin\cS$. This, according to Lemma \ref{lem:B}, implies $F(x)\cap N_\delta=\emptyset$. If $x\in N_\delta\setminus \ogr{\cS}$, then we again obtain $F(x)\cap N_\delta=\emptyset$ as a consequence of Lemma \ref{lem:C}. Both cases show that there is no solution with respect to $F$ passing through $x$ and contained in $N_\delta$, which means that $\Inv_FN_\delta\subset N_\epsilon\cap\ogr{\cS}$, and therefore completes the proof.
\qed

\medskip

Note that by Theorem \ref{thm:invNd=invNe} the sets $M_r$ can be alternatively expressed as
$$
M_r=\Inv_F(N^r _\delta),
$$
where
$
N^r _\delta=\bigcup_{\sigma\in\cM_r}\cl\csd{\sigma}.
$

\begin{lem}\label{lem:sol_prop}
Let $\phi:\ZZ\to X$ be a solution for the multivalued map $F$. Assume that the sequence of simplices $\rho^*:\ZZ\to\cX$ and $\rho:\ZZ\to\cX$ define a corresponding solution of the combinatorial flow~$\Pi_\cV$,
as introduced in Theorem~\ref{thm:orbitequiv2}.
If $\cS$ is an isolated invariant set with respect to $\Pi_\cV$, and if there exists a integer $k\in\ZZ$ such that $\rho^*(k), \rho^*(k+1)\in\cS$, and if each simplex in the extended solution $\rho$ between the simplices $\rho^*(k)$ and $\rho^*(k+1)$ belongs to $\cS$, then $\phi(k+1)\in\ogr{\cS}\cap N_\epsilon$.
\end{lem}
\proof
Observe that by \eqref{eq:orbitequiv2-1} we have the inclusion $\varphi(k+1)\in\cse{\rho^*(k+1)}$,
and since $\rho^*(k+1)\in\cS$  we get $\phi(k+1)\in N_\epsilon$.
We need to verify that $\phi (k+1)\in\ogr{\cS}$.
Let $\sigma_i=\sigma^\epsilon _{max}(\phi(i))=\rho^*(i)$  for $i\in\ZZ$.
Then we have to consider the following two complementary cases: $\phi(k+1)\in\sigma_k ^+$ and $\phi(k+1)\notin\sigma_k ^+$.

The first case immediately shows that $\phi(k+1)\in|\ccl\cS|$, as $\sigma_k\in\cS$ implies
the inclusion $\sigma_k ^+\in\cS$ according to the assertion that $\cS$ is an isolated invariant set (cf. Proposition~\ref{def:discrete-isolated-inv}). Since we also have $\phi(k+1)\in N_\epsilon$, the inclusion $\phi (k+1)\in\ogr{\cS}$ follows from Lemma \ref{lem:Ne_cap_supS_closed}.

Consider now the second case $\phi(k+1)\notin\sigma_k ^+$.
According to Theorem~\ref{thm:orbitequiv2}(1)
we have the inclusion $\phi(k+1)\in |\ccl\tau^+|$, where~$\tau$ is a simplex in the extended solution~$\rho$
which lies between the simplices $\rho^*(k)$ and $\rho^*(k+1)$. According to our assumption $\tau$ belongs to $\cS$, hence so does $\tau^+$, as $\cS$ is isolated and invariant. Consequently, $\phi(k+1)\in |\ccl\cS|$ which, along with the inclusion $\phi(k+1)\in N_\epsilon$ and Lemma \ref{lem:Ne_cap_supS_closed}, implies $\phi(k+1)\in \ogr{\cS}$. This completes the proof.\qed
\medskip

\subsection{Proof of Theorem \ref{thm:Mcomb=MF}.}
First note that the sets~$M_r$ of~$M$ are mutually disjoint, which is a consequence of the mutual disjointness of the sets in the family~$\cM$ and the definition of the sets $M_r$.
Moreover, by Theorem~\ref{th:N-iso-block} and Theorem \ref{thm:invNd=invNe}, they are isolated invariant sets with respect to the map~$F$. Hence, condition~(a) of Definition~\ref{def:Morse_dec} holds.

We now verify condition~(b) of Definition \ref{def:Morse_dec}. Let $\phi:\ZZ\to X$ be an arbitrary solution of the multivalued map $F$. Let $\rho:\ZZ\to\cX$ denote a corresponding solution of the multivalued flow~$\Pi_\cV$ as constructed in Theorem~\ref{thm:orbitequiv2}. Since~$\cM$ is a Morse decomposition of~$\cV$,
there exist two indices $r,r'\in\PP$ with $r' \geq r$, such that the inclusions $\alpha(\rho)\subset \cM_{r'}$ and $\omega(\rho)\subset \cM_r$ are satisfied.

Let us first focus on the $\omega$-limit set of $\phi$. The inclusion $\omega(\rho)\subset \cM_r$ implies that there exists a $k'\in\ZZ^+$ such that $\rho(n)\in\cM_{r}$ for all $n\geq k'$. Passing to the reduced solution $\rho^*$ we infer that $\rho^*(n)\in\cM_r$ for large enough $n$. Let $k\in\ZZ^+$ be such that both
$\rho(n)\in\cM_{r}$ and $\rho^*(n)\in\cM_{r}$ hold for $n\geq k$.
Then Lemma~\ref{lem:sol_prop} implies the inclusion $\phi(k+1)\in M_r$. Applying Lemma \ref{lem:sol_prop}
another time, we further obtain the inclusion $\phi([k+1,+\infty))\subset M_r$, which
in combination with the closedness of~$M_r$ yields $\omega(\phi)\subset M_r$.

For the set $\alpha(\phi)$ a similar argument applies. Indeed, the inclusion $\alpha(\rho)\subset \cM_{r'}$
establishes the existence of $k\in\ZZ^+$ with $\rho(-n)\in\cM_{r'}$ and $\rho^*(-n)\in\cM_{r'}$ for all
integers $n\geq k$. Applying Lemma \ref{lem:sol_prop}, this time to the arguments~$-(k+1)$ and~$-k$,
we can further deduce that $\phi(-k)\in M_{r'}$. Now, by the reverse recurrence and Lemma \ref{lem:sol_prop}, we have $\phi((-\infty,-k])\subset M_r$, and the inclusion $\alpha(\phi)\subset M_{r'}$ follows.

Next, we verify condition~(c) of Definition~\ref{def:Morse_dec}. Suppose that~$\phi$
is a full solution of~$F$ such that $\alpha(\phi)\cup\omega(\phi)\subset M_r$ for
some $r\in\PP$. Consider the corresponding solution $\rho:\ZZ\to\cX$ of the
multivalued flow $\Pi_\cV$, as constructed in Theorem~\ref{thm:orbitequiv2}.
Since $\cM$ is a Morse decomposition with respect to $\Pi_\cV$, we have
$\alpha(\rho)\subset \cM_{r_1}$ and $\omega(\rho)\subset \cM_{r_2}$, for
some $r_1,r_2\in\PP$. Then the argument used for the proof of condition~(b)
shows that the two inclusions $\alpha(\phi)\subset M_{r_1}$ and
$\omega(\phi)\subset M_{r_2}$ are satisfied. This immediately yields
$r_1=r_2=r$, as $M$ is a family of disjoint sets. Thus, we have
$\alpha(\rho)\cup\omega(\rho)\subset \cM_r$. Since $\cM$ is a Morse
decomposition, the inclusion $\im\rho\subset\cM_r$ follows, and consequently
$\rho(k)\in\cM_r$ and $\rho^*(k)\in\cM_r$ for all $k\in\ZZ$. Again, by the
recurrent argument with respect to $k$ in both forward and backward directions
and Lemma \ref{lem:sol_prop} we conclude that $\im\phi\subset M_r$. This
completes the proof that the collection $M$ is a Morse decomposition of~$X$
with respect to~$F$.

The Conley indices of $\cM_r$ and $M_r$ coincide by
Theorem~\ref{thm:ConleyF=ConleyV}. The fact that the Conley-Morse
graphs coincide as well follows from Theorems~\ref{thm:orbitequiv2}
and~\ref{thm:orbitequiv1}(a) via an argument similar to the argument
for condition~(b) of Definition~\ref{defn:redextsoln} and is left
to the reader.
\qed


\end{document}

%% file: mathDefs.tex
\def\refeq#1{\if\workingver y(\ref{#1})-[[#1]]\else(\ref{#1})\fi}
\def\refth#1{\if\workingver y\ref{#1}-[[#1]]\else\ref{#1}\fi}
\def\mylabel#1{\if\workingver y\label{#1}{\bf\ \ [[#1]]\ \ }\else\label{#1}\fi}
\def\mybibitem#1{\if\workingver y\bibitem{#1}{\bf\ \ [[#1]]\ \
}\else\bibitem{#1}\fi}

\newfont{\msam}{msam10}
\newfont{\msbm}{msbm10}

\def\articletheorems{
\newtheorem{thm}{Theorem}[section]
\newtheorem{lem}[thm]{Lemma}

\newtheorem{defn}[thm]{Definition}
\newtheorem{cor}[thm]{Corollary}
\newtheorem{prop}[thm]{Proposition}

\newtheorem{algo}[thm]{Algorithm}
\newtheorem{alg}[thm]{Algorithm}  

}

\def\omap{\oarrow}

\def\cA{\text{$\mathcal A$}}

\def\cF{\text{$\mathcal F$}}

\def\cM{\text{$\mathcal M$}}

\def\cP{\text{$\mathcal P$}}

\def\cS{\text{$\mathcal S$}}
\def\cT{\text{$\mathcal T$}}

\def\cV{\text{$\mathcal V$}}

\def\cX{\text{$\mathcal X$}}

\newcommand{\id}{\operatorname{id}}
\newcommand{\cl}{\operatorname{cl}}

\newcommand{\Int}{\operatorname{int}}
\newcommand{\inte}{\operatorname{int}}

\newcommand{\bd}{\operatorname{bd}}
\newcommand{\Bd}{\operatorname{Bd}}

\newcommand{\dom}{\operatorname{dom}}

\newcommand{\sgn}{\operatorname{sgn}}

\newcommand{\im}{\operatorname{im}}

\newcommand{\conv}{\protect\mbox{\rm conv\,}}

\newcommand{\Con}{\operatorname{Con}}
\newcommand{\Inv}{\operatorname{Inv}}
\newcommand{\Fix}{\operatorname{Fix}}

\def\proof{{\bf Proof:\ }}

\def\begeq#1{\begin{equation}\mylabel{#1}}
\def\endeq{\end{equation}}

\def\mathobj#1{\mbox{$#1$}}

\def\ZZ{\mathobj{\mathbb{Z}}}

\def\RR{\mathobj{\mathbb{R}}}

\def\implies{\;\Rightarrow\;}
\def\iff{\;\Leftrightarrow\;}

\def\setof#1{\mbox{$\{\,#1\,\}$}}

\def\0#1{\hbox{\kern25pt}$ #1 $\\}
\def\1#1{\hbox{\kern40pt}$ #1 $\\}
\def\2#1{\hbox{\kern55pt}$ #1 $\\}
\def\3#1{\hbox{\kern70pt}$ #1 $\\}

\newcounter{li}

\def\begalg#1{\begin{algo}\mylabel{#1}\normalshape:\small\baselineskip 10pt\\}
\def\endalg{\end{algo}}

\def\Figures(include=#1,cat=#2){
  \renewcommand{\textfraction}{.20}
  \renewcommand{\topfraction}{.80}
  \renewcommand{\bottomfraction}{.80}
  \renewcommand{\floatpagefraction}{.80}
  \newcount\figcount
  \figcount=0
  \let\includefigures=#1
  \def\figcat{#2}
}

\def\FigureFromFile[#1][#2](#3)#4
{
  \begin{figure}[htbp]
     \global\advance\figcount by 1
     \if\includefigures y\special{anisoscale #1.wmf, \the\hsize #2}\fi
     \vspace{#2}
     \caption{#4}
     \mylabel{#3}
   \end{figure}
}

\def\FigureFromFileTwoD[#1][#2,#3](#4)#5
{
  \begin{figure}[htbp]
     \global\advance\figcount by 1
     \if\includefigures y\special{anisoscale #1.wmf, #2 #3}\fi
     \vspace{#2}
     \caption{#5}
     \mylabel{#4}
   \end{figure}
}

\def\FigureF<#1>[#2](#3)#4
{
  \begin{figure}[htbp]
     \global\advance\figcount by 1
     \if\includefigures y\special{anisoscale \figcat/fig\number\figcount.wmf,
       \the\hsize #2}
     \fi
     \if\includefigures p
       \leavevmode
       \epsfxsize=\hsize
       \epsffile{#1}
     \fi
     \if\includefigures y
          \vspace{#2}
     \fi
     \caption{#4}
     \mylabel{#3}
   \end{figure}
}

\def\Figure[#1](#2)#3
{
  \begin{figure}[htbp]
     \global\advance\figcount by 1
     \if\includefigures y\special{anisoscale \figcat/fig\number\figcount.wmf,
       \the\hsize #1}
     \fi
     \if\includefigures p
       \leavevmode
       \epsfxsize=\hsize
       \epsffile{fig\number\figcount.eps}
     \fi
     \if\includefigures y
          \vspace{#1}
     \fi
     \caption{#3}
     \mylabel{#2}
   \end{figure}
}